\documentclass[12pt]{amsart}

\usepackage{amsmath}
\usepackage{amscd}
\usepackage{amssymb}
\usepackage{latexsym}
\usepackage{color}
\usepackage{cancel}
\usepackage[arrow,curve,matrix,tips,2cell]{xy}

\newtheorem{theorem}{Theorem}[section]
\newtheorem*{theorem*}{Theorem}
\newtheorem{lemma}[theorem]{Lemma}
\newtheorem*{lemma*}{Lemma}
\newtheorem{corollary}[theorem]{Corollary}
\newtheorem*{corollary*}{Corollary}
\newtheorem{proposition}[theorem]{Proposition}

\newtheorem{remark}[theorem]{Remark}
\newtheorem{question}[theorem]{Question}
\newtheorem{definition}[theorem]{Definition}

\newtheorem{example}[theorem]{Example}
\newcommand{\bgl}{\begin{equation}} 
\newcommand{\egl}{\end{equation}}
\newcommand{\bgloz}{\begin{equation*}} 
\newcommand{\egloz}{\end{equation*}}
\newcommand{\bgln}{\begin{eqnarray}} 
\newcommand{\egln}{\end{eqnarray}}
\newcommand{\bglnoz}{\begin{eqnarray*}} 
\newcommand{\eglnoz}{\end{eqnarray*}}
\newcommand{\btheo}{\begin{theorem}}
\newcommand{\etheo}{\end{theorem}}
\newcommand{\btheooz}{\begin{theorem*}}
\newcommand{\etheooz}{\end{theorem*}}
\newcommand{\blemma}{\begin{lemma}}
\newcommand{\elemma}{\end{lemma}}
\newcommand{\blemmaoz}{\begin{lemma*}}
\newcommand{\elemmaoz}{\end{lemma*}}
\newcommand{\bproof}{\begin{proof}}
\newcommand{\eproof}{\end{proof}}
\newcommand{\bbew}{\begin{beweis}}
\newcommand{\ebew}{\end{beweis}}
\newcommand{\bremark}{\begin{remark}\em}
\newcommand{\eremark}{\end{remark}}
\newcommand{\bquestion}{\begin{question}\em}
\newcommand{\equestion}{\end{question}}
\newcommand{\bdefin}{\begin{definition}}
\newcommand{\edefin}{\end{definition}}
\newcommand{\bprop}{\begin{proposition}}
\newcommand{\eprop}{\end{proposition}}
\newcommand{\bcor}{\begin{corollary}}
\newcommand{\ecor}{\end{corollary}}
\newcommand{\bcoroz}{\begin{corollary*}}
\newcommand{\ecoroz}{\end{corollary*}}
\newcommand{\bfa}{\begin{cases}} 
\newcommand{\efa}{\end{cases}}
\newcommand{\bexample}{\begin{example}\em}
\newcommand{\eexample}{\end{example}}
\newcommand{\cB}{\mathcal B}
\newcommand{\cC}{\mathcal C}

\newcommand{\cM}{\mathcal M}
\newcommand{\cN}{\mathcal N}

\newcommand{\cS}{\mathcal S}

\def\Cz{\mathbb{C}}
\def\Fz{\mathbb{F}}

\def\Nz{\mathbb{N}}

\def\Qz{\mathbb{Q}}
\def\Rz{\mathbb{R}}
\def\Tz{\mathbb{T}}
\def\Zz{\mathbb{Z}}

\def\1z{\mathbb{1}}

\newcommand{\an}[1]{``#1''} 
\newcommand{\ti}{\tilde}

\newcommand{\lori}{\longrightarrow}

\newcommand{\ma}{\mapsto} 
\newcommand\onto{\twoheadrightarrow} 
\newcommand\into{\hookrightarrow} 
\newcommand{\Rarr}{\Rightarrow} 
\newcommand{\Larr}{\Leftarrow} 
\newcommand{\LRarr}{\Leftrightarrow} 

\def\SEMI{\mbox{$\times\kern-2pt\vrule height5pt width.6pt \kern3pt $}}

\newcommand{\img}{{\rm Im\,}}

\newcommand{\id}{{\rm id}}

\renewcommand{\ker}{{\rm ker}\,}

\newcommand{\abs}[1]{\left|#1\right|} 
\newcommand{\defeq}{\mathrel{:=}} 

\newcommand{\dop}{\text{: }} 
\newcommand{\plim}{\varprojlim} 

\newcommand{\lge}{\left\{} 
\newcommand{\rge}{\right\}} 
\newcommand{\lru}{\left(} 
\newcommand{\rru}{\right)} 
\newcommand{\rukl}[1]{\lru #1 \rru} 
\newcommand{\gekl}[1]{\lge #1 \rge} 

\newcommand{\menge}[2]{\gekl{ #1 \dop #2 }} 
\title{Continuous orbit equivalence rigidity}

\author{Xin Li}

\address{Xin Li, School of Mathematical Sciences, Queen Mary University of London, Mile End Road, London E1 4NS, United Kingdom}
\email{xin.li@qmul.ac.uk}

\subjclass[2010]{Primary 37B05; Secondary 37A20, 46L05}

\begin{document}

\begin{abstract}
We take the first steps towards a better understanding of continuous orbit equivalence, i.e., topological orbit equivalence with continuous cocycles. First, we characterise continuous orbit equivalence in terms of isomorphisms of C*-crossed products preserving Cartan subalgebras. This is the topological analogue of the classical result by Singer and Feldman-Moore in the measurable setting. Secondly, we turn to continuous orbit equivalence rigidity, i.e., the question whether for certain classes of topological dynamical systems, continuous orbit equivalence implies conjugacy. We show that this is not always the case by constructing topological dynamical systems (actions of free abelian groups, and also non-abelian free groups) which are continuously orbit equivalent but not conjugate. Furthermore, we prove positive rigidity results. For instance, it turns out that general topological Bernoulli actions are rigid when compared with actions of nilpotent groups, and that topological Bernoulli actions of duality groups are rigid when compared with actions of solvable groups. The same is true for certain subshifts of full shifts over finite alphabets.
\end{abstract}

\thanks{Research supported by EPSRC grant EP/M009718/1.}

\maketitle


\setlength{\parindent}{0cm} \setlength{\parskip}{0.5cm}

\section{Introduction}

From its very beginning on, the theory of operator algebras was closely related to ergodic theory and dynamical systems. The bridge between these subjects is built by crossed product constructions, attaching von Neumann algebras to measure-preserving dynamical systems and C*-algebras to topological dynamical systems.

In the setting of von Neumann algebras, the crossed product construction, also called group-measure space construction, played an important role in the classification of injective factors. Similarly, in the C*-algebraic setting, crossed products attached to topological dynamical systems provide interesting examples which are challenging to classify and lead to new insights.

If we want to further develop the relationship between operator algebras and dynamical systems, the following question will be crucial:
\setlength{\parindent}{0.5cm} \setlength{\parskip}{0cm}

How much information do these crossed product constructions contain about the underlying dynamical systems?
\setlength{\parindent}{0cm} \setlength{\parskip}{0cm}

It turns out that the crossed product itself might contain very little information, but if we consider the crossed product together with a commutative subalgebra (which is canonically given), then our question can be answered in a systematic way.
\setlength{\parindent}{0cm} \setlength{\parskip}{0.5cm}

To explain this in the measurable and von Neumann algebraic setting, let $G \curvearrowright X$ and $H \curvearrowright Y$ be probability measure preserving actions. Here, our measure spaces are standard, our groups are discrete and countable, and they act by Borel automorphisms. We say that $G \curvearrowright X$ and $H \curvearrowright Y$ are orbit equivalent if there exists an isomorphism of measure spaces $\varphi: \: X \to Y$ with $\varphi(G.x) = H.\varphi(x)$ for a.e. $x \in X$. Moreover, we let $G \ltimes X$ and $H \ltimes Y$ be the measured transformation groupoids attached to $G \curvearrowright X$ and $H \curvearrowright Y$. If our actions are (essentially) free, then $G \ltimes X$ and $H \ltimes Y$ are nothing else but the orbit equivalence relations $R(G \curvearrowright X)$ and $R(H \curvearrowright Y)$ viewed as measured groupoids. Here is a classical result:
\btheo[\cite{Sin, FMI, FMII}]
\label{vN}
Let $G \curvearrowright X$ and $H \curvearrowright Y$ be (essentially) free probability measure preserving actions. The following are equivalent:
\begin{itemize}
\item $G \curvearrowright X$ and $H \curvearrowright Y$ are orbit equivalent;
\item $G \ltimes X$ and $H \ltimes Y$ are isomorphic as measured groupoids (or equivalence relations);
\item there is a vN-isomorphism $\Phi: \: L^{\infty}(X) \rtimes G \overset{\cong}{\lori} L^{\infty}(Y) \rtimes H$ with $\Phi(L^{\infty}(X)) = L^{\infty}(Y)$.
\end{itemize}
\etheo
The interested reader may consult \cite{Vaes07, Ga, Vaes10} for more details.

Our first result carries over Theorem~\ref{vN} to the topological setting. Let $G \curvearrowright X$ and $H \curvearrowright Y$ be topological dynamical systems. This means that $G$ and $H$ are countable discrete groups acting by homeomorphisms on locally compact Hausdorff spaces $X$ and $Y$.
\setlength{\parindent}{0.5cm} \setlength{\parskip}{0cm}

We say that $G \curvearrowright X$ and $H \curvearrowright Y$ are continuously orbit equivalent if there exists a homeomorphism $\varphi: \: X \overset{\cong}{\lori} Y$ with inverse $\psi = \varphi^{-1}: \: Y \overset{\cong}{\lori} X$ and continuous maps $a: G \times X \to H$, $b: H \times Y \to G$ such that $\varphi(g.x) = a(g,x).\varphi(x)$ and $\psi(h.y) = b(h,y).\psi(y)$ for all $g \in G$, $x \in X$, $h \in H$ and $y \in Y$. Note that $G$ and $H$ carry the discrete topology. This notion of continuous orbit equivalence has been studied in special cases (see \cite{GPS, BT, Su}), but not -- at least to the best of the author's knowledge -- in the general setting. There is also a weaker notion of (topological) orbit equivalence which has been studied intensively for $\Zz^n$-actions on the Cantor set in the remarkable papers \cite{GPS, GMPS08, GMPS10}.

Moreover, let $G \ltimes X$ and $H \ltimes Y$ be the transformation groupoids attached to $G \curvearrowright X$ and $H \curvearrowright Y$. Here is the topological analogue of Theorem~\ref{vN}:
\setlength{\parindent}{0cm} \setlength{\parskip}{0.5cm}
\btheo
\label{COE-gpd-Cartan}
Let $G \curvearrowright X$ and $H \curvearrowright Y$ be topologically free systems, and assume that $X$ and $Y$ are second countable. The following are equivalent:
\begin{itemize}
\item $G \curvearrowright X$ and $H \curvearrowright Y$ are continuously orbit equivalent;
\item $G \ltimes X$ and $H \ltimes Y$ are isomorphic as topological groupoids;
\item there is a C*-isomorphism $\Phi: \: C_0(X) \rtimes_r G \overset{\cong}{\lori} C_0(Y) \rtimes_r H$ with $\Phi(C_0(X)) = C_0(Y)$.
\end{itemize}
\etheo
Here and in the sequel, \an{topologically free system} stands for \an{topologically free topological dynamical system}. Because of Theorem~\ref{COE-gpd-Cartan}, it seems that -- at least for our purposes -- continuous orbit equivalence is a good topological analogue of orbit equivalence in the measurable setting.

In the measurable setting, orbit equivalence rigidity has established itself as a key notion. The idea is to find classes of actions for which orbit equivalence already implies conjugacy. Indeed, impressive orbit equivalence rigidity results have been obtained in \cite{Zim, Fur, MS, Popa07, Popa08, Kida, Ioana}. Viewing continuous orbit equivalence as the topological analogue of orbit equivalence, a natural question is whether there are rigidity phenomena for continuous orbit equivalence.

The only result known in this context is due to \cite{BT}, which says that if $\Zz \curvearrowright X$ and $\Zz \curvearrowright Y$ are topologically free systems on compact spaces $X$ and $Y$ such that one of them is topologically transitive, then $\Zz \curvearrowright X$ and $\Zz \curvearrowright Y$ must already be conjugate if they are continuously orbit equivalent. Apart from this, not much else seems to be known about continuous orbit equivalence rigidity.

The main goal of the present paper is to take the first steps towards a better understanding of continuous orbit equivalence rigidity.

First of all, we construct examples of topological dynamical systems which are continuously orbit equivalent but not conjugate. This ensures that the comparison between continuous orbit equivalence and conjugacy is really interesting. A first class of examples is given by products of odometer actions. A second family of examples is constructed from boundary actions of non-abelian free groups and odometer actions, inspired by \cite{Su}. 

Secondly, we prove positive results in continuous orbit equivalence rigidity.
\btheo
\label{main1}
Let $G$ be a torsion-free group, and let $H$ be a finitely generated nilpotent group which is not virtually infinite cyclic. Assume that $G \curvearrowright X$ is a topologically free system on a compact space $X$ such that $G \curvearrowright X$ is almost $\Zz G$-projective. Furthermore, let $H \curvearrowright Y$ be a topologically free system. If $G \curvearrowright X$ and $H \curvearrowright Y$ are continuously orbit equivalent, then they must be conjugate. 
\etheo
\btheo
\label{main2}
Let $G$ be a duality group in the sense of \cite[Chapter~VIII, \S~10]{Bro} which is not infinite cyclic, and let $H$ be a finitely generated solvable group. Assume that $G \curvearrowright X$ is a topologically free system on a compact space $X$ such that $G \curvearrowright X$ is almost $\Zz G$-projective. Furthermore, let $H \curvearrowright Y$ be a topologically free system. If $G \curvearrowright X$ and $H \curvearrowright Y$ are continuously orbit equivalent, then they must be conjugate. 
\etheo
Here, we say that $G \curvearrowright X$ is almost $\Zz G$-projective if $C(X,\Zz) \cong \Zz \oplus P$ as $\Zz G$-modules, where the copy of $\Zz$ is given by the constant functions on $X$ and $P$ is a projective $\Zz G$-module. For instance, as we will see, the Bernoulli action $G \curvearrowright X_0^G$ is almost $\Zz G$-projective for every compact space $X_0$ and every torsion-free group $G$. Also, for a torsion-free group, a subshift of the full shift over a finite alphabet whose forbidden words avoid a fixed letter is almost $\Zz G$-projective. These systems are actually even almost $\Zz G$-free, in the sense that $P$ can be chosen to be $\Zz G$-free. This leads to the following immediate consequence:
\bcor
Let $G \curvearrowright X$ be a topological Bernoulli action with compact base or a subshift of the full shift over a finite alphabet whose forbidden words avoid a fixed letter. Moreover, let $H \curvearrowright Y$ be a topologically free system. Assume that $G$ is a torsion-free group and that $H$ is a finitely generated nilpotent group which is not virtually infinite cyclic, or that $G$ is a duality group which is not infinite cyclic and that $H$ is a finitely generated solvable group. If $G \curvearrowright X$ and $H \curvearrowright Y$ are continuously orbit equivalent, then they must be conjugate.
\ecor

In view of these results, an interesting and natural task is to find more examples of topological dynamical systems which are almost $\Zz G$-projective. Apart from Bernoulli actions and subshifts, it turns out that Denjoy homeomorphisms, restricted to their unique minimal closed invariant subspaces, give rise to Cantor minimal systems which are almost $\Zz G$-free (where $G = \Zz$).

Theorem~\ref{main1} and Theorem~\ref{main2} are our main results. Their proofs consist of three main ingredients. Each of them is interesting in its own right.

The first ingredient establishes a link between continuous cocycle rigidity and continuous orbit equivalence rigidity. Let $G \curvearrowright X$ be a topological dynamical system, and let $H$ be a group. A continuous function $a: \: G \times X \to H$ is called a continuous $H$-cocycle for $G \curvearrowright X$ if $a(g_1 g_2,x) = a(g_1, g_2.x) a(g_2,x)$ for all $g_1, g_2 \in G$ and $x \in X$. In particular, we can view any group homomorphism $\rho: \: G \to H$ as a cocycle given by $(g,x) \ma \rho(g)$. Continuous cocycles $a$ and $a'$ are called cohomologous if there exists a continuous map $u: X \to H$ such that $a(g,x) = u(g.x) a'(g,x) u(x)^{-1}$. We say that $G \curvearrowright X$ is continuous $H$-cocycle rigid if every continuous $H$-cocycle for $G \curvearrowright X$ is cohomologous to some group homomorphism $\rho: \: G \to H$.
\setlength{\parindent}{0cm} \setlength{\parskip}{0cm}

In general, it is not clear how continuous cocycle rigidity is related to continuous orbit equivalence rigidity. However, we have
\btheo
\label{CCR->COER}
Let $G$ be a torsion-free amenable group. Assume that $G \curvearrowright X$ and $H \curvearrowright Y$ are topologically free systems on compact spaces $X$ and $Y$, and suppose that $G \curvearrowright X$ and $H \curvearrowright Y$ are continuously orbit equivalent. If $G \curvearrowright X$ is continuous $H$-cocycle rigid, then $G \curvearrowright X$ and $H \curvearrowright Y$ must be conjugate.
\etheo
\setlength{\parindent}{0cm} \setlength{\parskip}{0.5cm}

The second ingredient establishes continuous cocycle rigidity for certain actions and certain groups.
\btheo
\label{CCR}
Let $G$ be a duality group in the sense of \cite[Chapter~VIII, \S~10]{Bro} with $cd(G) \neq 1$, let $X$ be a compact space, and suppose that $G \curvearrowright X$ is a topological dynamical system which is almost $\Zz G$-projective. Then $G \curvearrowright X$ is continuous $H$-cocycle rigid for every solvable group $H$.
\etheo

The third and final ingredient builds a bridge between continuous orbit equivalence and the notion of quasi-isometry or rather topological couplings. Inspiration comes from \cite{Fur, Sha, Sau}.
\setlength{\parindent}{0cm} \setlength{\parskip}{0cm}

Given groups $G$ and $H$, a topological coupling $G \curvearrowright \Omega \curvearrowleft H$ for $G$ and $H$ is a locally compact Hausdorff space $\Omega$ with a left $G$-action commuting with a right $H$-action such that both actions admit compact open fundamental domains. A topological coupling $G \curvearrowright \Omega \curvearrowleft H$ is called topologically free if the corresponding $G \times H$-action is topologically free.

A continuous orbit couple for $G$ and $H$ consists of topological dynamical systems $G \curvearrowright X$, $H \curvearrowright Y$ on compact spaces $X$, $Y$, continuous maps $p: \: X \to Y$, $q: \: Y \to X$, $a: \: G \times X \to H$, $b: \: H \times Y \to G$, $g: \: X \to G$ and $h: \: Y \to H$ such that $p(g.x) = a(g,x).p(x)$, $q(h.y) = b(h,y).q(y)$, $qp(x) = g(x).x$ and $pq(y) = h(y).y$. We call a continuous orbit couple for $G$ and $H$ topologically free if $G \curvearrowright X$ and $H \curvearrowright Y$ are topologically free. If $G \curvearrowright X$ and $H \curvearrowright Y$ are continuously orbit equivalent via $\varphi:\: X \cong Y$, then we obtain a continuous orbit couple for $G$ and $H$ by setting $p = \varphi$ and $q = \varphi^{-1}$. 
\btheo
\label{1-1}
Let $G$ and $H$ be groups. Then there exists a one-to-one correspondence between (isomorphism classes of) topologically free continuous orbit couples for $G$ and $H$ and (isomorphism classes of) topologically free topological couplings for $G$ and $H$.
\etheo
\setlength{\parindent}{0cm} \setlength{\parskip}{0.5cm}

In particular, if topological dynamical systems $G \curvearrowright X$, $H \curvearrowright Y$ on compact spaces $X$, $Y$ are continuously orbit equivalent, and if one of the groups ($G$ or $H$) is finitely generated, then $G$ and $H$ must be quasi-isometric.

In \S~\ref{sec-COE-gpd-Cartan}, we introduce the notion of continuous orbit equivalence, make some general observations and prove Theorem~\ref{COE-gpd-Cartan}. Moreover, we discuss known examples for continuous orbit equivalence rigidity and construct counterexamples for which continuous orbit equivalence does not imply conjugacy in \S~\ref{sec-ex}. In \S~\ref{sec-CCR->COER}, we introduce the notion of continuous cocycle rigidity, study the connection to continuous orbit equivalence rigidity, and prove Theorem~\ref{CCR->COER}. In the following section (\S~\ref{sec-CCR-cohom}), we show that Bernoulli actions, certain subshifts as well as Cantor minimal systems arising from Denjoy homeomorphisms are almost $\Zz G$-free, study continuous cocycle rigidity using non-abelian group cohomology and prove Theorem~\ref{CCR}. Thereafter, we introduce the notions of topological couplings and continuous orbit couples, establish the connection between them, and prove Theorem~\ref{1-1}. Finally, in \S~\ref{sec-concl}, we prove Theorem~\ref{main1} and Theorem~\ref{main2}.

I would like to thank David Kerr for inspiring discussions about continuous orbit equivalence.

\section{Continuous orbit equivalence, transformation groupoids and Cartan pairs}
\label{sec-COE-gpd-Cartan}

All our groups are discrete and countable, and all our topological spaces are locally compact and Hausdorff. By a topological dynamical system, we mean an action of a group on a topological space by homeomorphisms. In this section, all our topological spaces are second countable.

Let $G \curvearrowright X$ be a topological dynamical system. The $G$-action is denoted by $G \times X \to X, \, (g,x) \ma g.x$. For $x \in X$, let $G_x = \menge{g \in G}{g.x = x}$ be its stabilizer group. The transformation groupoid $G \ltimes X$ attached to $G \curvearrowright X$ is given by the set $G \times X$ with multiplication $(g',x')(g,x) = (g'g,x)$ if $x' = g.x$, inversion $(g,x)^{-1} = (g^{-1},g.x)$, range map $r(g,x) = g.x$ and source map $s(g,x) = x$. Obviously, $G \ltimes X$ is {\'e}tale. The reduced groupoid C*-algebra $C^*_r(G \ltimes X)$ is canonically isomorphic to $C_0(X) \rtimes_r G$. Moreover, we have a canonical embedding $C_0(X) \into C_0(X) \rtimes_r G$.

\bdefin
$G \curvearrowright X$ is called topologically free if for every $e \neq g \in G$, $\menge{x \in X}{g.x \neq x}$ is dense in $X$.
\edefin
From now on, for the sake of brevity, we write \an{topologically free system} for \an{topologically free topological dynamical system}.

\blemma
\label{TF--dense}
$G \curvearrowright X$ is topologically free if and only if $\menge{x \in X}{G_x = \gekl{e}}$ is dense in $X$.
\elemma
\bproof
\an{$\Larr$} is clear. For \an{$\Rarr$}, note that by topological freeness, $\menge{x \in X}{g.x \neq x}$ is dense (and open) in $X$ for all $e \neq g \in G$. Thus, by the Baire category theorem,
$$
  \menge{x \in X}{G_x = \gekl{e}} = \bigcap_{e \neq g \in G} \menge{x \in X}{g.x \neq x}
$$
must be dense in $X$.
\eproof

\bcor
\label{TF--tp}
$G \curvearrowright X$ is topologically free if and only if the transformation groupoid $G \ltimes X$ is topologically principal.
\ecor
\bproof
By definition (see \cite{R08}), $G \ltimes X$ is topologically principal if and only if the set of points in $X$ with trivial isotropy is dense in $X$. But this set coincides with $\menge{x \in X}{G_x = \gekl{e}}$. Thus Lemma~\ref{TF--dense} implies our corollary.
\eproof

\bremark
Corollary~\ref{TF--tp} shows that if $G \curvearrowright X$ is topologically free, then the pair $(C_0(X) \rtimes_r G, C_0(X))$ is a Cartan pair in the sense of \cite[Definition~5.1]{R08}.
\eremark

Recall the following definition from the introduction:
\bdefin
\label{defCOE}
Topological dynamical systems $G \curvearrowright X$ and $H \curvearrowright Y$ are continuously orbit equivalent (we write $G \curvearrowright X \ \sim_{\text{\tiny coe}} \ H \curvearrowright Y$) if there exists a homeomorphism $\varphi: \: X \overset{\cong}{\lori} Y$ with inverse $\psi = \varphi^{-1}: \: Y \overset{\cong}{\lori} X$ and continuous maps $a: G \times X \to H$, $b: H \times Y \to G$ such that
\bgln
\label{phigx}
  \varphi(g.x) &=& a(g,x).\varphi(x) \\
\label{psihy}
  \psi(h.y) &=& b(h,y).\psi(y)
\egln
for all $g \in G$, $x \in X$, $h \in H$ and $y \in Y$.
\edefin

\bremark
\label{GxHy}
\eqref{phigx} implies $\varphi(G.x) \subseteq H.\varphi(x)$ for all $x \in X$, and \eqref{psihy} implies $\psi(H.y) \subseteq G.\psi(y)$ for all $y \in Y$. Thus, $\varphi(G.x) = H.\varphi(x)$ and $\psi(H.y) = G.\psi(y)$.
\eremark

\bremark
If $H \curvearrowright Y$ is topologically free, then $a$ is uniquely determined by \eqref{phigx}, and by symmetry, if $G \curvearrowright X$ is topologically free, then $b$ is uniquely determined by \eqref{psihy}. The reason is as follows: Suppose that $a': \: G \times X \to H$ is another continuous map with $\varphi(g.x) = a'(g,x).\varphi(x)$. For arbitrary $g \in G$ and $x \in X$, there exists an open neighbourhood $U$ of $x$ such that $a$ and $a'$ are constant on $\gekl{g} \times U$, with values $h$ and $h'$ in $H$, say. Then for every $\bar{x} \in U$, $\varphi(g.\bar{x}) = h.\varphi(\bar{x}) = h'.\varphi(\bar{x})$. Topological freeness implies $h = h'$, in particular $a(g,x) = a'(g,x)$.
\eremark

\blemma
\label{cocycle-id}
In Definition~\ref{defCOE}, if $H \curvearrowright Y$ is topologically free, then 
$$a(g_1g_2,x) = a(g_1,g_2.x) a(g_2,x)$$
for all $g_1, g_2 \in G$ and $x \in X$.
\elemma
\bproof
Let $g_1, g_2 \in G$ and $x \in X$ be arbitrary. Choose an open neighbourhood $U$ of $x \in X$ such that $a(g_1g_2,\bar{x}) = a(g_1g_2,x)$, $a(g_1,g_2.\bar{x}) = a(g_1,g_2.x)$ and $a(g_2,\bar{x}) = a(g_2,x)$ for all $\bar{x} \in U$. Then for all $\bar{x} \in U$, $\varphi(g_1g_2.\bar{x}) = \varphi(g_1.(g_2.\bar{x})) = a(g_1,g_2.\bar{x}).\varphi(g_2.\bar{x}) = a(g_1,g_2.\bar{x}) a(g_2,\bar{x}).\varphi(\bar{x}) = a(g_1,g_2.x) a(g_2,x).\varphi(\bar{x})$, but also $\varphi(g_1g_2.\bar{x}) = a(g_1g_2,\bar{x}).\varphi(\bar{x}) = a(g_1g_2,x).\varphi(\bar{x})$. By topological freeness, $a(g_1g_2,x) = a(g_1,g_2.x) a(g_2,x)$.
\eproof

\blemma
\label{a_x-bij}
In the situation of Definition~\ref{defCOE}, let $Y_f = \menge{y \in Y}{H_y = \gekl{e}}$. For every $x \in \psi(Y_f)$, $a_x: \: G \to H, \, g \ma a(g,x)$ is bijective.
\elemma
\bproof
Since $\varphi(x) \in Y_f$, $a_x$ is injective. To prove surjectivity, take $h \in H$. Since by Remark~\ref{GxHy}, $\varphi(G.x) = H.\varphi(x)$, there exists $g \in G$ with $h.\varphi(x) = \varphi(g.x) = a(g,x).\varphi(x)$. As $\varphi(x) \in Y_f$, we conclude that $h = a(g,x) = a_x(g)$.
\eproof

\blemma
\label{lem-baphi}
In the situation of Definition~\ref{defCOE}, assume that $G \curvearrowright X$ and $H \curvearrowright Y$ are topologically free. Then
\bgl
\label{baphi}
  b(a(g,x),\varphi(x)) = g \ {\rm for} \ {\rm all} \ g \in G, \, x \in X,
\egl
and $b$ is uniquely determined by \eqref{baphi}.
\elemma
\bproof
Let $h \defeq a(g,x)$. Then $\varphi(g.x) = h.\varphi(x)$, so $g.x = \psi(h.\varphi(x)) = b(h,\varphi(x)).x$. Since this equation holds in an open neighbourhood of $x$, topological freeness implies $b(a(g,x),\varphi(x)) = g$. Moreover, note that for all $x \in \psi(Y_f)$, $a_x(G) = H$ by Lemma~\ref{a_x-bij}. Hence \eqref{baphi} determines $b$ on $H \times Y_f$. But since $Y_f$ is dense in $Y$ by topological freeness, and because $b$ is continuous, \eqref{baphi} determines $b$ on $H \times Y$.
\eproof

\bcor
\label{XfYf}
In the situation of Definition~\ref{defCOE}, assume that $G \curvearrowright X$ and $H \curvearrowright Y$ are topologically free. Let $X_f = \menge{x \in X}{G_x = \gekl{e}}$ and $Y_f = \menge{y \in Y}{H_y = \gekl{e}}$. Then $\varphi(X_f) = Y_f$. In particular, for every $x \in X$ with $G_x = \gekl{e}$, $a_x: \: G \to H, \, g \ma a(g,x)$ is bijective.
\ecor
\bproof
By symmetry, we just have to show $\varphi(X_f) \subseteq Y_f$. Take $x \in X_f$, and let $y = \varphi(x)$. Suppose that $h \in H$ satisfies $h.y = y$. Then $x = \psi(y) = \psi(h.y) = b(h,y).\psi(y) = b(h,y).x$, and therefore $b(h,y) = e$ since $x \in X_f$. But by the analogue of \eqref{baphi} with reversed roles for $a$ and $b$, we get $e = a(e,x) = a(b(h,y),x) = h$. Hence $y \in Y_f$.
\eproof

We are now ready for the proof of Theorem~\ref{COE-gpd-Cartan}.
\btheooz[Theorem~\ref{COE-gpd-Cartan}]
Let $G \curvearrowright X$ and $H \curvearrowright Y$ be topologically free systems. The following are equivalent:
\begin{itemize}
\item[(i)] $G \curvearrowright X \ \sim_{\text{\tiny coe}} \ H \curvearrowright Y$;
\item[(ii)] $G \ltimes X \cong H \ltimes Y$ (as topological groupoids);
\item[(iii)] there is a C*-isomorphism $\Phi: \: C_0(X) \rtimes_r G \overset{\cong}{\lori} C_0(Y) \rtimes_r H$ with $\Phi(C_0(X)) = C_0(Y)$.
\end{itemize}
\etheooz
\bproof
(i) $\Rarr$ (ii): Assume that $G \curvearrowright X \ \sim_{\text{\tiny coe}} \ H \curvearrowright Y$, and let $\varphi$, $\psi$, $a$ and $b$ be as in Definition~\ref{defCOE}. Then $G \ltimes X \to H \ltimes Y, \, (g,x) \ma (a(g,x),\varphi(x))$ and $H \ltimes Y \to G \ltimes X, \, (h,y) \ma (b(h,y),\psi(y))$ are certainly continuous groupoid morphisms, and they are inverse to each other due to \eqref{baphi} and the analogue of \eqref{baphi} with reversed roles for $a$ and $b$.

(ii) $\Rarr$ (i): Let $\chi: \: G \ltimes X \overset{\cong}{\lori} H \ltimes Y$ be an isomorphism of topological groupoids. Set $\varphi = \chi \vert_X: \: X \overset{\cong}{\lori} Y$ and let $a$ be the composition $G \ltimes X \overset{\chi}{\lori} H \ltimes Y \to H$, where the second map is $H \ltimes Y \to H, (h,y) \ma h$. Then $a$ is obviously continuous, and $\varphi(g.x) = \chi(r(g,x)) = r(\chi(g,x)) = r(a(g,x),\varphi(x)) = a(g,x).\varphi(x)$. Similarly, for $\psi = \varphi^{-1}$, if we let $b$ be the composition $H \ltimes Y \overset{\chi^{-1}}{\lori} G \ltimes X \to G$, where the second map is $G \ltimes X \to G, (g,x) \ma g$, then $\psi(h.y) = b(h,y).\psi(y)$.

(ii) $\LRarr$ (iii) is \cite[Proposition~4.13]{R08}, where we have to use Corollary~\ref{TF--tp}.
\eproof

\section{Continuous orbit equivalence rigidity: Examples and counterexamples}
\label{sec-ex}

Let us compare continuous orbit equivalence with conjugacy.
\bdefin
Topological dynamical systems $G \curvearrowright X$ and $H \curvearrowright Y$ are conjugate (we write $G \curvearrowright X \ \sim_{\text{\tiny conj}} \ H \curvearrowright Y$) if there is a homeomorphism $\varphi: \: X \overset{\cong}{\lori} Y$ and a group isomorphism $\rho: \: G \overset{\cong}{\lori} H$ such that for every $g \in G$ and $x \in X$, $\varphi(g.x) = \rho(g).\varphi(x)$.
\edefin

Obviously, $G \curvearrowright X \ \sim_{\text{\tiny conj}} \ H \curvearrowright Y$ implies $G \curvearrowright X \ \sim_{\text{\tiny coe}} \ H \curvearrowright Y$. Are there classes of dynamical systems where we can reverse this implication, i.e., where continuous orbit equivalence implies conjugacy?

Here is a first class of examples, for which continuous orbit equivalence rigidity holds because of a trivial reason: Suppose that $G \curvearrowright X$ is a topologically free system on a connected space $X$. If $G \curvearrowright X \ \sim_{\text{\tiny coe}} \ H \curvearrowright Y$ for some topologically free system $H \curvearrowright Y$, then $G \curvearrowright X \ \sim_{\text{\tiny conj}} \ H \curvearrowright Y$. The reason is that the function $a$ in Definition~\ref{defCOE} is continuous, hence for every $g \in G$, $a \vert_{\gekl{g} \times X}$ is constant because $X$ is connected and $H$ is discrete. Hence $a(g,x) = \rho(g)$ for some map $\rho: \: G \to H$, and $\rho$ has to be a homomorphism (by Lemma~\ref{cocycle-id}) and bijective (by Lemmma~\ref{a_x-bij}).

This observation means that if we focus on discrete groups, it is natural to restrict our discussion to topological dynamical systems on totally disconnected spaces.

Here is a first result in continuous orbit equivalence rigidity:
\btheo[{\cite[Theorem~3.2]{BT}}]
Let $\Zz \curvearrowright X$ and $\Zz \curvearrowright Y$ be topologically free systems on compact spaces $X$ and $Y$. Assume that $\Zz \curvearrowright X$ is topologically transitive.

If $\Zz \curvearrowright X \ \sim_{\text{\tiny coe}} \ \Zz \curvearrowright Y$, then $\Zz \curvearrowright X \ \sim_{\text{\tiny conj}} \ \Zz \curvearrowright Y$.
\etheo
In this theorem, while the groups are fixed, the assumptions on the actions are very mild. Therefore, an immediate question is whether there are counterexamples to continuous orbit equivalence rigidity at all, i.e., examples of topological dynamical systems which are continuously orbit equivalent but not conjugate.

\subsection{Products of odometer transformations}

Let $M = \prod_p p^{v_p}$ be a supernatural number. Here, the product is taken over all primes, $v_p \in \gekl{0, 1, 2, \dotsc} \cup \gekl{\infty}$, and $\sum_p v_p = \infty$. The odometer action $\Zz \curvearrowright \Zz / M$ corresponding to $M$ is constructed as follows:

Choose a sequence $(m_k)_k$ of natural numbers such that, for all primes $p$, $v_p(m_k) \nearrow v_p$ for $k \to \infty$. Then set $\Zz / M = \plim_k \Zz / m_k$. The canonical projections $\Zz \onto \Zz / m_k$ induce a group embedding $\Zz \into \Zz / M$, and this in turn yields an action $\Zz \curvearrowright \Zz / M$ which we call the odometer transformation for $M$.

\btheo
\label{odo-coe}
For supernatural numbers $M_1, \dotsc, M_r$ and $N_1, \dotsc, N_s$, the following are equivalent:
\begin{enumerate}
\item[(i)] $\Zz^r \curvearrowright \prod_{i=1}^r \Zz / M_i \ \sim_{\text{\tiny coe}} \ \Zz^s \curvearrowright \prod_{j=1}^s \Zz / N_j$;
\item[(ii)] $C_0(\prod_{i=1}^r \Zz / M_i) \rtimes \Zz^r \cong C_0(\prod_{j=1}^s \Zz / N_j) \rtimes \Zz^s$;
\item[(iii)] $(K_*(C_0(\prod_{i=1}^r \Zz / M_i) \rtimes \Zz^r), [1]_0) \cong (K_*(C_0(\prod_{j=1}^s \Zz / N_j) \rtimes \Zz^s), [1]_0)$;
\item[(iv)] $r=s$, there exists $\sigma \in S_r$, natural numbers $m_1, \dotsc, m_r$ and $n_1, \dotsc, n_r$ such that for all $1 \leq i \leq r$, $m_i M_i = n_{\sigma(i)} N_{\sigma(i)}$, and $\prod_{i=1}^r M_i = \prod_{j=1}^r N_j$.
\end{enumerate}
\etheo
\bproof
(i) $\Rarr$ (ii) follows from Theorem~\ref{COE-gpd-Cartan} as our systems are free.

(ii) $\Rarr$ (iii) is clear.

(iii) $\Rarr$ (iv): $K_*$ stands for $K_0 \oplus K_1$. Clearly, $(K_0(C(\Zz / M) \rtimes \Zz),[1]_0) \cong (\Zz[M^{-1}],1)$ and $K_1(C(\Zz / M) \rtimes \Zz) \cong \Zz$. Here $\Zz[M^{-1}] = \menge{\frac{x}{m} \in \Qz}{m \mid M}$. So $K_*(C(\prod_{i=1}^r \Zz / M_i) \rtimes \Zz^r) \cong \bigoplus_{I \subseteq \gekl{1, \dotsc, r}} \Zz[(\prod_{i \in I} M_i)^{-1}]$ and $[1]_0$ corresponds to $1 \in \Zz[(\prod_{i=1}^r M_i)^{-1}]$ ($I = \gekl{1, \dotsc, r}$).

Therefore, $\Qz^{2^r} \cong K_*(C(\prod_{i=1}^r \Zz / M_i) \rtimes \Zz^r) \otimes \Qz \cong K_*(C(\prod_{j=1}^s \Zz / N_j) \rtimes \Zz^s) \otimes \Qz \cong \Qz^{2^s}$, and this implies $r=s$. Moreover, as a $K_*$-isomorphism preserves $[1]_0$, it restricts to an isomorphism $\Zz[(\prod_{i=1}^r M_i)^{-1}] \cong \Zz[(\prod_{j=1}^r N_j)^{-1}]$ sending $1$ to $1$. This implies $\prod_{i=1}^r M_i = \prod_{j=1}^r N_j$.

Given supernatural numbers $M$ and $N$, we define $M \lesssim N$ if there exists $n \in \Nz$ with $M \mid nN$ ($v_p(M) \leq v_p(nN)$). We define $M \sim N$ if $M \lesssim N$ and $N \lesssim M$. It is immediate that there exists a non-zero homomorphism $\Zz[M^{-1}] \to \Zz[N^{-1}]$ if and only if $M \lesssim N$. Set $\cM = \menge{M_i}{1 \leq i \leq r}$, $\bigwedge \cM = \menge{\prod_{i \in I} M_i}{I \subseteq \gekl{1,\dotsc,r}}$ and $\cN = \menge{N_j}{1 \leq j \leq r}$, $\bigwedge \cN = \menge{\prod_{j \in J} N_j}{J \subseteq \gekl{1,\dotsc,r}}$. Using the assumption that $\bigoplus_{I \subseteq \gekl{1, \dotsc, r}} \Zz[(\prod_{i \in I} M_i)^{-1}] \cong \bigoplus_{J \subseteq \gekl{1, \dotsc, r}} \Zz[(\prod_{j \in J} N_j)^{-1}]$, a straightforward inductive argument shows that for every equivalence class $\cS$ of supernatural numbers with respect to $\sim$, $\abs{\cS \cap \bigwedge \cM} = \abs{\cS \cap \bigwedge \cN}$, and then also $\abs{\cS \cap \cM} = \abs{\cS \cap \cN}$.

(iv) $\Rarr$ (i): We need the following observation: Let $l$ be a natural number and $\lambda_l: \: \Zz / l \curvearrowright \Zz / l$ the canonical action. Let $L$ be a supernatural number, $X = \Zz / lL$, $\ti{X} = l \cdot (\Zz / lL)$, $\alpha_{lL}: \: \Zz \curvearrowright X$ the odometer transformation for $lL$, and $\ti{\alpha} = \alpha \vert_{L \Zz}: l \Zz \curvearrowright \ti{X}$. We claim that
\bgl
\label{alphalambdaalpha}
\alpha_{lL} \ \sim_{\text{\tiny coe}} \ \lambda_l \boxtimes \ti{\alpha}_{lL} \ \sim_{\text{\tiny conj}} \ \lambda_l \boxtimes \alpha_L.
\egl
$\boxtimes$ stand for the product action. Let us prove \eqref{alphalambdaalpha}. Define $\varphi: \: X = \bigsqcup_{k=0}^{l-1} k + \ti{X} \to \Zz / l \times \ti{X}, \, k + x \ma ([k],x)$. It is easy to see that the inverse of $\varphi$ is given by $\psi: \: \Zz / l \times \ti{X} \to X, \, ([k],x) \ma k + x$ for $0 \leq k \leq l-1$. Moreover, define $a: \: \Zz \times X = \bigsqcup_{j=0}^{l-1} (j+l\Zz) \times \bigsqcup_{k=0}^{l-1} (k + \ti{X}) \to \Zz / l \times l \Zz$ by setting $a(j+h,k+x) = ([j],h)$ if $j+k \leq l-1$ and $a(j+h,k+x) = ([j],k+l)$ if $l < j+k$. Also, define $b: \: (\Zz / l \times l \Zz) \times (\Zz / l \times \ti{X}) \to \Zz$ by setting $b(([j],h),([k],x)) = j+h$ if $j+k \leq l-1$ and $b(([j],h),([k],x)) = j+h-l$ if $l \leq j+k$, where $0 \leq j, k \leq l-1$. Then it is easy to check that $\varphi$, $a$, $\psi$ and $b$ satisfy \eqref{phigx} and \eqref{psihy}, so that $\alpha_{lL} \ \sim_{\text{\tiny coe}} \ \lambda_l \boxtimes \ti{\alpha}_{lL}$. Furthermore, $\lambda_l \boxtimes \ti{\alpha}_{lL} \ \sim_{\text{\tiny conj}} \ \lambda_l \boxtimes \alpha_L$ is easy to see. This proves \eqref{alphalambdaalpha}.

Now we can complete the proof for (iv) $\Rarr$ (i). Without loss of generality we may assume that $\sigma = \id$, i.e., $m_i M_i = n_i N_i$ for all $1 \leq i \leq r$. Without loss of generality, we may further assume $\gcd(m_i,n_i) = 1$. Then we can write $M_i = n_i L_i$ and $N_i = m_i L_i$ for some supernatural number $L_i$. Set $L = \prod_{i=1}^r L_i$, and choose natural numbers $m$ and $n$ with $\gcd(m,L) = 1 = \gcd(n,L)$ such that $\prod_{i=1}^r M_i = (\prod_{i=1}^r n_i)(\prod_{i=1}^r L_i) = nL$ and $\prod_{j=1}^r N_j = (\prod_{j=1}^r m_j)(\prod_{j=1}^r L_j) = mL$. $\prod_{i=1}^r M_i = \prod_{j=1}^r N_j$ implies that $m=n$. Therefore, we get
\bglnoz
  \boxtimes_{i=1}^r \alpha_{M_i} &=& \boxtimes_{i=1}^r \alpha_{n_i L_i} \ \overset{\eqref{alphalambdaalpha}}{\sim}_{\text{\tiny coe}} \ \boxtimes_{i=1}^r (\lambda_{n_i} \boxtimes \alpha_{L_i}) \ \overset{\eqref{alphalambdaalpha}}{\sim}_{\text{\tiny coe}} \ \alpha_{n L_1} \boxtimes \alpha_{L_2} \boxtimes \dotso \boxtimes \alpha_{L_r} \\
  &=& \alpha_{m L_1} \boxtimes \alpha_{L_2} \boxtimes \dotso \boxtimes \alpha_{L_r} \overset{\eqref{alphalambdaalpha}}{\sim}_{\text{\tiny coe}} \boxtimes_{j=1}^r \alpha_{N_j}.
\eglnoz
\eproof

In contrast, for conjugacy, we get
\btheo
\label{odo-conj}
Let $I$ and $J$ be finite sets. For supernatural numbers $\gekl{M_i}_{i \in I}$ and $\gekl{N_j}_{j \in J}$, $\Zz^I \curvearrowright \prod_{i \in I} \Zz / M_i \ \sim_{\text{\tiny conj}} \ \Zz^J \curvearrowright \prod_{j \in J} \Zz / N_j$ if and only if there exists a finite set $K$, supernatural numbers $\gekl{L_k}_{k \in K}$ such that
\begin{itemize}
\item $I = \bigsqcup_{k \in K} I_k$, $J = \bigsqcup_{k \in K} J_k$,
\item $\abs{I_k} = \abs{J_k}$ for all $k \in K$,
\item for every $k \in K$, every $i \in I_k$ and $j \in J_k$, we can write $M_i = m_i L_k$ and $N_j = n_j L_k$ for some (uniquely determined) $m_i, n_j \in \Nz$ with $\gcd(m_i,L_k) = 1 = \gcd(n_j,L_k)$, and we have $\prod_{i \in I_k} \Zz / m_i \cong \prod_{j \in J_k} \Zz / n_j$.
\end{itemize}
\etheo
\bproof
\an{$\Rarr$}: Assume that $\rho: \: \Zz^I \cong \Zz^J$ is a group isomorphism and $\varphi: \: \prod_{i \in I} \Zz / M_i \cong \prod_{j \in J} \Zz / N_j$ such that $\varphi(g.x) = \rho(g).\varphi(x)$. $\Zz^I \cong \Zz^J$ implies that $\abs{I} = \abs{J}$. Set $r = \abs{I} = \abs{J}$. Moreover, we may assume $\varphi(0) = 0$ (otherwise go over to $\varphi - \varphi(0)$). Let $\rho$ be multiplication with $S \in GL_r(\Zz)$.

It is straightforward to check that if $S_{j,i} \neq 0$, then $N_j \lesssim M_i$. So there exist a finite set $K$ and decompositions $I = \bigsqcup_{k \in K} I_k$, $J = \bigsqcup_{k \in K} J_k$ with $\abs{I_k} = \abs{J_k}$ for all $k \in K$ such that for every $(i,j) \in I_k \times J_k$, $M_i \sim N_j$.

Fix $k \in K$. Find a supernatural number $L_k$ such that for every $i \in I_k$, $j \in J_k$, there are $m_i \in \Nz$, $n_j \in \Nz$ with $\gcd(m_i,L_k) = 1 = \gcd(n_j,L_k)$ such that $M_i = m_i L_k$ and $N_j = n_j L_k$. Then $\varphi$ restricts to an isomorphism of topological abelian groups
$$
  \varphi_k: (\prod_{i \in I_k} \Zz / m_i) \times (\prod_{i \in I_k} \Zz / L_k) \cong \prod_{i \in I_k} \Zz / M_i \overset{\cong}{\lori} \prod_{j \in J_k} \Zz / N_j \cong (\prod_{j \in J_k} \Zz / n_j) \times (\prod_{j \in J_k} \Zz / L_k).
$$
Let $l \in \Nz$ satisfy $\prod_{j \in J_k} n_j \mid l$ and $\gcd(l,L) = 1$. Certainly, $\varphi_k(w,0)$ is of the form $(z,0)$ as for all $i \in I_k$, there exists no non-zero homomorphism $\Zz / m_i \to \Zz / L_k$. Also, $\varphi_k(0,y)$ is of the form $\varphi_k(0,l \cdot \ti{y}) = l \cdot \varphi_k(0,\ti{y})$, hence of the form $(0,x)$ as for all $j \in J_k$, $l \equiv 0$ in $\Zz / n_j$. Hence $\varphi_k = \phi_t \times \phi_L$ for some group isomorphisms $\phi_t: \: \prod_{i \in I_k} \Zz / m_i \overset{\cong}{\lori} \prod_{j \in J_k} \Zz / n_j$ and $\phi_L: \: \prod_{i \in I_k} \Zz / L_k \overset{\cong}{\lori} \prod_{j \in J_k} \Zz / L_k$.

\an{$\Larr$}: Without loss of generality we may assume $\abs{K} = 1$. Let $K = \gekl{k}$, $I = I_k$, $J = J_k$, $\abs{I} = \abs{J} = r$. We may assume that $I = J = \gekl{1, \dotsc, r}$. Let $L = L_k$ be a supernatural number such that for every $1 \leq i, j \leq r$, $M_i = m_i L$ and $N_j = n_j L$ for some (unique) $m_i, n_j \in \Nz$ with $\gcd(m_i,L) = 1 = \gcd(n_j,L)$, and such that $\prod_{i=1}^r \Zz / m_i \cong \prod_{j=1}^r \Zz / n_j$. By the theory of elementary divisors, there are $S, T \in GL_r(\Zz)$ such that
$
S
\rukl{
\begin{smallmatrix}
m_1& & 0 \\ & \ddots & \\ 0 & & m_r
\end{smallmatrix}
}
T
=
\rukl{
\begin{smallmatrix}
n_1& & 0 \\ & \ddots & \\ 0 & & n_r
\end{smallmatrix}
}
$. Thus
$
S
\rukl{
\rukl{
\begin{smallmatrix}
m_1& & 0 \\ & \ddots & \\ 0 & & m_r
\end{smallmatrix}
}
\Zz^r
}
=
\rukl{
\begin{smallmatrix}
n_1& & 0 \\ & \ddots & \\ 0 & & n_r
\end{smallmatrix}
}
\Zz^r
$.
So the same matrix $S$ induces two group isomorphisms $\rho: \: \Zz^r \to \Zz^r, \, g \ma Sg$ and 
$
\phi_t: \: \prod_{i=1}^r \Zz / m_i \cong
\Zz^r {\bigg /}
\rukl{
\begin{smallmatrix}
m_1& & 0 \\ & \ddots & \\ 0 & & m_r
\end{smallmatrix}
} 
\Zz^r
\cong
\Zz^r {\bigg /}
\rukl{
\begin{smallmatrix}
n_1& & 0 \\ & \ddots & \\ 0 & & n_r
\end{smallmatrix}
} 
\Zz^r
\cong \prod_{i=j}^r \Zz / n_j
$
as well as an isomorphism of topological groups $\phi_L: \: (\Zz / L)^r \cong (\Zz / L)^r$. Let $\varphi$ be the isomorphism
\bglnoz
  && \prod_{i=1}^r \Zz / M_i \cong \prod_{i=1}^r (\Zz / m_i \times \Zz / L) \cong (\prod_{i=1}^r \Zz / m_i) \times (\Zz / L)^r \\
  &\overset{\phi_t \times \phi_L}{\lori}&
  (\prod_{j=1}^r \Zz / n_j) \times (\Zz / L)^r \cong \prod_{j=1}^r (\Zz / n_j \times \Zz / L) \cong \prod_{j=1}^r \Zz / N_j.
\eglnoz
Then $\varphi \vert_{\Zz^r} = \rho$, and so $\varphi(g.x) = \rho(g).\varphi(x)$ for all $g \in \Zz^r$ and $x \in \prod_{i=1}^r \Zz / M_i$. This means that $\Zz^I \curvearrowright \prod_{i \in I} \Zz / M_i \ \sim_{\text{\tiny conj}} \ \Zz^J \curvearrowright \prod_{j \in J} \Zz / N_j$.
\eproof

Comparing Theorem~\ref{odo-coe} and Theorem~\ref{odo-conj}, we can easily construct products of odometers which are continuously orbit equivalent but not conjugate.
\bexample
Let $r \geq 2$. Let $p$ and $q$ be primes, $p \neq q$, and let $n \in \Nz$ with $n > 1$ and $\gcd(p,n) = 1 = \gcd(q,n)$. If we set $M_1 = n \cdot p^{\infty}$, $M_2 = q^{\infty}$, $M_3 = \dotso = M_s = p^{\infty}$ and $N_1 = p^{\infty}$, $N_2 = n \cdot q^{\infty}$, $N_3 = \dotso = N_s = p^{\infty}$, then $\Zz^r \curvearrowright \prod_{i=1}^r \Zz / M_i \ \sim_{\text{\tiny coe}} \ \Zz^r \curvearrowright \prod_{j=1}^r \Zz / N_j$ but $\Zz^r \curvearrowright \prod_{i=1}^r \Zz / M_i \ \cancel{\sim}_{\text{\tiny conj}} \ \Zz^r \curvearrowright \prod_{j=1}^r \Zz / N_j$.
\eexample

\subsection{Actions of non-abelian free groups}

Let us construct actions of the free group $\Fz_r$ ($r \geq 2$) on the Cantor set, which are continuously orbit equivalent but not conjugate. Let $a_1, \dotsc, a_r$ be generators of $\Fz_r$. Let $\beta: \: \Fz_r \curvearrowright \partial \Fz_r$ be the $\Fz_r$-action on the Gromov boundary of $\Fz_r$, and set $\beta_i \defeq \beta_{a_i}$. For a supernatural number $M$, let $\alpha_M: \: \Zz \curvearrowright \Zz / M$ be the corresponding odometer transformation. For supernatural numbers $M_1$, $M_2$, $N_1$ and $N_2$, define actions $\gamma: \: \Fz_r \curvearrowright \partial \Fz_r \times (\Zz / M_1) \times (\Zz / M_2)$ and $\delta: \: \Fz_r \curvearrowright \partial \Fz_r \times (\Zz / N_1) \times (\Zz / N_2)$ by setting $\gamma_1 \defeq \beta_1 \times \alpha_{M_1} \times \id$, $\gamma_2 \defeq \beta_2 \times \id \times \alpha_{M_2}$, $\gamma_i \defeq \beta_i \times \id \times \id$ for all $i \geq 3$, $\gamma_{a_i} = \gamma_i$ for all $1 \leq i \leq r$, and similarly $\delta_1 \defeq \beta_1 \times \alpha_{N_1} \times \id$, $\delta_2 \defeq \beta_2 \times \id \times \alpha_{N_2}$, $\delta_i \defeq \beta_i \times \id \times \id$ for all $i \geq 3$, $\delta_{a_i} = \delta_i$ for all $1 \leq i \leq r$.
\btheo
\label{F-action}
Let $p$ and $q$ be primes, $p \neq q$, and let $n \in \Nz$ with $n > 1$ and $\gcd(p,n) = 1 = \gcd(q,n)$. If we set $M_1 = n \cdot p^{\infty}$, $M_2 = q^{\infty}$ and $N_1 = p^{\infty}$, $N_2 = n \cdot q^{\infty}$, then $\gamma \ \sim_{\text{\tiny coe}} \ \delta$ but $\gamma \ \cancel{\sim}_{\text{\tiny conj}} \ \delta$.
\etheo
For the proof, we need some preparation. Let $X$ be a totally disconnected compact space. In our application, $X$ will be the Cantor space. Let $C^{\infty}(X,\Cz) = C(X,\Zz) \otimes \Cz$. Obviously, we have an isomorphism $C^{\infty}(X,\Cz) \overset{\cong}{\lori} \menge{X \overset{f}{\lori} \Cz \ {\rm continuous}}{f(X) \subseteq \Cz \ {\rm finite}}, \, f \otimes z \ma f \cdot z$. Here $(f \cdot z)(x) = f(x) z$. In the following, we view elements in $C^{\infty}(X,\Cz)$ as $\Cz$-valued continuous functions on $X$ via this explicit isomorphism. Let $\phi: \: X \to X$ be a homeomorphism, and denote the induced automorphism of $C(X)$ by $\phi$ again. Obviously, $\phi(C^{\infty}(X,\Cz)) \subseteq C^{\infty}(X,\Cz)$. We define $E(\phi) \defeq \menge{z \in \Tz}{\phi(f) = z f \ {\rm for} \ {\rm some} \ 0 \neq f \in C^{\infty}(X,\Cz)}$. Now let $g_1, \dotsc, g_r$ be generators of $\Fz_r$, and let $g = g_1$. Let $Y$ be a totally disconnected compact space and let $\alpha: \: Y \cong Y$ be a homeomorphism.
\bprop
\label{EphiEalpha}
For $\phi = \beta_g \times \alpha$, we have $E_{\phi} = E_{\alpha}$.
\eprop
\bproof
We think of elements in $\partial \Fz_r$ as infinite reduced words in $g^{\pm 1} = g_1^{\pm 1}, g_2^{\pm 1}, \dotsc, g_r^{\pm 1}$. Let $W$ be the set of finite reduced words in $g_1^{\pm 1}, g_2^{\pm 1}, \dotsc, g_r^{\pm 1}$ which do not end on $g_r^{-1}$ nor on $g_r^2$. For $w \in W$, let $C_w$ be the subspace of $\partial \Fz_r$ consisting of those infinite reduced words which start with $w$. Note that the empty word $\emptyset$ lies in $W$, and that $C_{\emptyset} = \partial \Fz_r$. Clearly, $\menge{1_{C_w}}{w \in W}$ is a $\Zz$-basis for $C(\partial \Fz_r,\Zz)$. Take a family $\cC$ of compact open subsets of $Y$ such that $\menge{1_C}{C \in \cC}$ is a $\Zz$-basis for $C(Y,\Zz)$. Then $\menge{1_{C_w} \otimes 1_C}{w \in W, C \in \cC}$ is a $\Zz$-basis for $C(\partial \Fz_r \times Y,\Zz) \cong C(\partial \Fz_r,\Zz) \otimes C(Y,\Zz)$.

Let $z \in E_{\phi}$, and let $f = \sum_i 1_{C_{w_i}} \otimes 1_{C_i} \otimes \lambda_i$ ($\lambda_i \in \Cz \setminus \gekl{0}$) be a non-zero element in $C(\partial \Fz_r \times Y,\Cz) \cong C(\partial \Fz_r,\Zz) \otimes C(Y,\Zz) \otimes \Cz$ with $\phi(f) = zf$. Then $\sum_i 1_{C_{w_i}} \otimes (1_{C_i} \otimes z \lambda_i) = zf = \phi(f) = \sum_i 1_{\beta_g(C_{w_i})} \otimes (1_{\alpha(C_i)} \otimes \lambda_i) = \sum_j 1_{C_{w_j}} \otimes f_j$ for some $0 \neq f_j \in C(Y,\Zz) \otimes \Cz$, where $\gekl{w_j}_j$ are chosen so that $1_{\beta_g(C_{w_i})} \in \Zz \text{-}{\rm span}(\gekl{w_j}_j)$ for all $i$. It follows that $\gekl{w_i}_i = \gekl{w_j}_j \supseteq \menge{g w_i}{w_i \neq g^{-1}} \cup \gekl{g}$ if $g^{-1} \in \gekl{w_i}_i$ and $\gekl{w_i}_i = \gekl{w_j}_j \supseteq \menge{g w_i}{w_i \neq g^{-1}}$ if $g^{-1} \notin \gekl{w_i}_i$. Here we use that $\beta_g(C_{w_i}) = C_{g w_i} \neq C_{w_i}$ if $w_i \neq g^{-1}$ and $\beta_g(C_{g^{-1}}) = \partial \Fz_r \setminus C_g$.

We claim that it already follows that $\gekl{w_i}_i = \gekl{\partial \Fz_r}$: If there is $w \in \gekl{w_i}_i$ not starting with $g^{-1}$, then $g^n w \in \gekl{w_i}$ for all $m \in \Nz$ which is impossible since $\gekl{w_i}_i$ is finite. If there is $w \in \gekl{w_i}_i$ of the form $g^{-m}v$ where $v \neq \emptyset$ is a finite reduced word not starting with $g^{\pm 1}$, then $v \in \gekl{w_i}_i$ contradicting our first observation. If there is $g^{-m} \in \gekl{w_i}_i$ for some $m \geq 1$, then $g^{-1} \in \gekl{w_i}_i$, hence $g \in \gekl{w_i}_i$. This again contradicts our first observation. Therefore, the only possibility is $\gekl{w_i}_i = \gekl{\partial \Fz_r}$.

Hence $f = 1 \otimes \ti{f}$ for some $\ti{f} \in C(Y,\Zz) \otimes \Cz$. So $1 \otimes z \ti{f} = zf = \phi(f) = 1 \otimes \alpha(\ti{f})$. Hence it follows that $z \ti{f} = \alpha(\ti{f})$. This shows that $z \in E_{\alpha}$. Since $z \in E_{\phi}$ was arbitrary, we obtain $E_{\phi} \subseteq E_{\alpha}$. The reverse inclusion is obvious.
\eproof
We are now ready for the
\bproof[Proof of Theorem~\ref{F-action}]
By \cite[Theorem~5.10 and Proposition~5.2]{Su}, $\gamma \ \sim_{\text{\tiny coe}} \ \delta$. So we just have to show $\gamma \ \cancel{\sim}_{\text{\tiny conj}} \ \delta$. Assume that there exists a homeomorphism $\varphi: \: \partial \Fz_r \times (\Zz / M_1) \times (\Zz / M_2) \overset{\cong}{\lori} \partial \Fz_r \times (\Zz / N_1) \times (\Zz / N_2)$ and a group isomorphism $\rho: \Fz_r \cong \Fz_r$ such that $\varphi \circ \gamma_a = \delta_{\rho(a)} \circ \varphi$ for all $a \in \Fz_r$. Let $a \defeq a_1$ and $g \defeq \rho(a_1)$. Then in particular, $\varphi \circ \gamma_a = \delta_g \circ \varphi$, so that $\gamma_a$ and $\delta_g$ are conjugate, and hence $E(\gamma_a) = E(\delta_g)$.

By construction, there are $k, l \in \Zz$ with $\delta_g = \beta_g \times \alpha_{N_1}^k \times \alpha_{N_2}^l$, where $N_1 = p^{\infty}$ and $N_2 = n \cdot q^{\infty}$. Proposition~\ref{EphiEalpha} yields $E(\delta_g) = E(\alpha_{N_1}^k \times \alpha_{N_2}^l)$. For a supernatural number $M$, let $\Tz(M) = \Zz[M^{-1}] / \Zz \subseteq \Qz / \Zz \subseteq \Rz / \Zz \cong \Tz$. If $l = 0$, then $E(\delta_g) = E(\alpha_{N_1}^k \times \id) = E(\alpha_{N_1}^k) = E(\alpha_{p^{\infty}}^k)$ is equal to $\gekl{1}$ or $\Tz(p^{\infty})$. If $l \neq 0$, then $\Tz(q^{\infty}) = E(\alpha_{q^{\infty}}^l) \subseteq E(\alpha_{n \cdot q^{\infty}}^l) \subseteq E(\alpha_{p^{\infty}}^k \times \alpha_{n \cdot q^{\infty}}^l) = E(\alpha_{N_1}^k \times \alpha_{N_2}^l) = E(\delta_g)$. However, $E(\gamma_a) = E(\alpha_{M_1} \times \id) = E(\alpha_{n \cdot p^{\infty}}) = \Tz(n \cdot p^{\infty})$ is not equal to $\gekl{1}$ nor $\Tz(p^{\infty})$ and does not contain $\Tz(q^{\infty})$. Hence $E(\gamma_a) \neq E(\delta_g)$. This is a contradiction.
\eproof

\section{From continuous cocycle rigidity to continuous orbit equivalence rigidity}
\label{sec-CCR->COER}

We introduce the notion of continuous cocycle rigidity. Let $G \curvearrowright X$ be a topological dynamical system and let $H$ be a group.
\bdefin
\label{def-cHc}
A continuous $H$-cocycle for $G \curvearrowright X$ is a continuous map $a: \: G \times X \to H$ such that $a(g_1 g_2,x) = a(g_1,g_2.x) a(g_2,x)$ for all $g_1, g_2 \in G$, $x \in X$.
\edefin
In other words, $a: \: G \times X \to H$ is a groupoid homomorphism, where we view $G \times X$ as a groupoid by identifying it with the transformation groupoid $G \ltimes X$ attached to $G \curvearrowright X$, and view $H$ as the groupoid whose unit space is a point.
\bdefin
\label{def-cohom}
Continuous $H$-cocycles $a$ and $a'$ for $G \curvearrowright X$ are continuously cohomologous ($a \sim a'$) if there exists a continuous map $u: \: X \to H$ such that $a(g,x) = u(g.x) a'(g,x) u(x)^{-1}$ for all $g \in G$ and $x \in X$.
\edefin
\bdefin
$G \curvearrowright X$ is continuous $H$-cocycle rigid if for every continuous $H$-cocycle $a$ for $G \curvearrowright X$, there exists a group homomorphism $\rho: \: G \to H$ such that $a \sim \rho$.
\edefin
Here we view $\rho$ as the cocycle $G \times X \to H, \, (g,x) \ma \rho(g)$.

The following observation provides a first link between continuous cocycle rigidity and continuous orbit equivalence rigidity.
\bprop
\label{prop_CCR->COER}
Let $G \curvearrowright X$ and $H \curvearrowright Y$ be topologically free systems. Assume that $G \curvearrowright X \ \sim_{\text{\tiny coe}} \ H \curvearrowright Y$, and let $\varphi$, $\psi$, $a$ and $b$ be as in Definition~\ref{defCOE}. If there exists a continuous map $u: \: X \to H$ and a group isomorphism $\rho: \: G \to H$ such that $a(g,x) = u(g.x) \rho(g) u(x)^{-1}$ for all $g \in G$, $x \in X$, then ${}_u \varphi: \: X \to Y, \, x \ma u(x)^{-1}.\varphi(x)$ and $\rho$ give rise to a conjugacy between $G \curvearrowright X$ and $H \curvearrowright Y$.
\eprop
\bproof
${}_u \varphi$ is obviously continuous, and an easy computation shows that ${}_u \varphi(g.x) = \rho(g).{}_u \varphi(x)$ for all $g \in G$, $x \in X$. It remains to show that ${}_u \varphi$ is a homeomorphism.

Let $\sigma = \rho^{-1}$, define $v: \: Y \to G, \, y \ma \rukl{\sigma(u(\psi(y)))}^{-1}$ and $\ti{b}: \: H \times Y \to G, \, (h,y) \ma v(h.y) \sigma(h) v(y)^{-1}$. Since
\bglnoz
  && \ti{b}(a(g,x),\varphi(x))
  = v(a(g,x).\varphi(x)) \sigma(a(g,x)) v(\varphi(x))^{-1}
  = v(\varphi(g.x)) \sigma(a(g,x)) v(\varphi(x))^{-1} \\
  &=& \sigma(u(g.x))^{-1} \sigma(a(g,x)) \sigma(u(x))
  = \sigma(u(g.x))^{-1} a(g,x) u(x))
  = \sigma(\rho(g)) = g,
\eglnoz
Lemma~\ref{lem-baphi} implies that $b = \ti{b}$. Set ${}_v \psi: \: Y \to X, \, y \ma v(y)^{-1} \psi(y)$. ${}_v \psi$ is obviously continuous, and an easy computation shows that ${}_v \psi(h.y) = \rho(h).{}_v \psi(y)$ for all $h \in H$, $y \in Y$. Moreover,
\bglnoz
  {}_v \psi({}_u \varphi(x)) &=& {}_v \psi(u(x)^{-1}.\varphi(x)) 
  = \sigma(u(x)^{-1}).{}_v \psi(\varphi(x)) \\
  &=& \sigma(u(x))^{-1} v(\varphi(x))^{-1}.x = \sigma(u(x))^{-1} \sigma(u(x)).x = x, \\
  and \
  {}_u \varphi({}_v \psi(y)) &=& {}_u \varphi(v(y)^{-1}.\psi(y)) 
  = \rho(v(y)^{-1}) {}_u \varphi(\psi(y)) \\
  &=& \rho(v(y))^{-1} u((\psi(y))^{-1}.y = u(\psi(y)) u((\psi(y))^{-1}.y = y.
\eglnoz
Thus ${}_u \varphi$ is a homeomorphism, with inverse ${}_v \psi$, and the proof is complete.
\eproof

Continuous cocycle rigidity means that every cocycle, whether or not if comes from a continuous orbit equivalence, is continuously cohomologous to a group homomorphism. At the same time, the preceding proposition shows that for continuous orbit equivalence rigidity, cocycles are required to be continuously cohomologous to group isomorphisms. Therefore, there does not seem to be any obvious connection between continuous cocycle rigidity and continuous orbit equivalence rigidity. However, we have the following
\btheo
\label{coe-rho-isom}
Suppose that $G$ is amenable and torsion-free. Let $G \curvearrowright X$ and $H \curvearrowright Y$ be topologically free systems on compact spaces $X$ and $Y$. Assume that $G \curvearrowright X \ \sim_{\text{\tiny coe}} \ H \curvearrowright Y$, and let $a: \: G \times X \to H$ be as in Definition~\ref{defCOE}. If $a \sim \rho$ for some group homomorphism $\rho: \: G \to H$, then $\rho$ must be an isomorphism.
\etheo
\bproof
Let $u: \: X \to H$ be continuous such that $a(g,x) = u(g.x) \rho(g) u(x)^{-1}$ for all $g \in G$, $x \in X$. Take $x \in X$ with $G_x = \gekl{e}$. Then by Lemma~\ref{XfYf}, $a_x: \: G \to H, \, g \ma a(g,x)$ is bijective. Let $u_x: \: G \to H, \, g \ma u(g.x)$. Then $a_x(g) = u_x(g) \rho(g) u_x(e)^{-1}$. $u$ is continuous and $X$ is compact, hence $u(X) \subseteq H$ is finite. In particular, $u_x(G)$ is finite. Therefore, for every $g \in \ker(\rho)$, $a_x(g) \in u_x(G) u_x(e)^{-1}$. It follows that $a_x(\ker(\rho))$ is finite. Since $a_x$ is injective, $\ker(\rho)$ is finite. But $G$ is torsion-free. This implies $\ker(\rho) = \gekl{e}$, so that $\rho$ is injective.

It remains to prove surjectivity. Since $a_x$ is surjective, we have $H = u(X) \rho(G) u(x)^{-1} = u(X) \rho(G)$. Thus, $[H:\rho(G)] < \infty$. In particular, $H$ is also amenable. Without loss of generality, we may assume $u_x(e) = u(x) = e$. Otherwise, replace $\rho$ by $u(x) \rho u(x)^{-1}$ and $u_x$ by $u_x \cdot u(x)^{-1}$. Suppose that $\rho(G) \subsetneq H$. Let $R$ be a complete system of left coset representatives of $\rho(G)$ in $H$. Since $H$ is amenable, there exists a finite subset $F$ of $H$ such that $\abs{rF \triangle F} < \frac{1}{3} \abs{F}$ for all $r \in R$ and $\abs{sF \triangle F} < \frac{1}{3 \abs{u(X)}} \abs{F}$ for all $s \in u(X)$.

Assume that $\abs{F \cap \rho(G)} > \frac{2}{3} \abs{F}$. By assumption ($\rho(G) \subsetneq H$), there exists $r \in R$ with $r \rho(G) \cap \rho(G) = \emptyset$. So $r (F \cap \rho(G)) \cap (F \cap \rho(G)) = \emptyset$, and we obtain $\abs{r(F \cap \rho(G)) \cap F} < \frac{1}{3} \abs{F}$. Moreover, $\abs{rF \setminus r(F \cap \rho(G))} = \abs{F \setminus (F \cap \rho(G))} < \frac{1}{3} \abs{F}$. Therefore, $\frac{2}{3} \abs{F} < \abs{rF \cap F} = \abs{(rF \setminus r(F \cap \rho(G))) \cup (r(F \cap \rho(G)) \cap F)} < \frac{1}{3} \abs{F} + \frac{1}{3} \abs{F} = \frac{2}{3} \abs{F}$. But this is a contradiction. Therefore, we must have $\abs{F \cap \rho(G)} \leq \frac{2}{3} \abs{F}$.

We certainly have $a_x^{-1}(F) \subseteq \rho^{-1}(\bigcup_{s \in u(X)} s^{-1}F)$. Hence 
\bglnoz
  \abs{F} &=& \abs{a_x^{-1}(F)} \leq \abs{\rho^{-1}(\bigcup_{s \in u(X)} s^{-1}F)} 
  \leq \abs{\rho^{-1}(F \cup (\bigcup_{s \in u(X)} s^{-1}F) \setminus F)} \\
  &\leq& \abs{\rho^{-1}(F)} + \sum_{s \in u(X)} \abs{s^{-1}F \setminus F} < \abs{F \cap \rho(G)} + \frac{1}{3} \abs{F} \leq \abs{F}.
\eglnoz
This is a contradiction. We conclude that $\rho(G) = H$.
\eproof

Clearly, Proposition~\ref{prop_CCR->COER} and Theorem~\ref{coe-rho-isom} imply
\btheooz[Theorem~\ref{CCR->COER}]
Let $G$ be a torsion-free amenable group. Assume that $G \curvearrowright X$ and $H \curvearrowright Y$ are topologically free systems on compact spaces $X$ and $Y$, and suppose that $G \curvearrowright X$ and $H \curvearrowright Y$ are continuously orbit equivalent. If $G \curvearrowright X$ is continuous $H$-cocycle rigid, then $G \curvearrowright X$ and $H \curvearrowright Y$ must be conjugate.
\etheooz

\section{Continuous cocycle rigidity via group cohomology}
\label{sec-CCR-cohom}

The first goal of this section is to rephrase continuous cocycle rigidity in the language of non-abelian group cohomology. For the sake of completeness, we briefly recall the definition of non-abelian group cohomology ($H^1$). We refer the reader to \cite[Part Three, Appendix \an{Non-abelian cohomology}]{Se1979} and \cite[Chapter~I, \S~5]{Se1997} for details.

Let $G$ be a group acting on a group $A$ by automorphisms, denoted by $G \times A \to A, \, (s,a) \ma s.a$. A $1$-cocycle of $G$ in $A$ is a map $G \to A, \, s \ma a_s$ such that $a_{st} = a_s s.a_t$. We write $Z^1(G,A)$ for the set of all these $1$-cocycles. Given $1$-cocycles $a$ and $a'$ of $G$ in $A$, we say that $a$ is cohomologous to $a'$ ($a \sim a'$) if there exists $b \in A$ with $a'_s = b^{-1} a_s s.b$ for all $s \in G$. We define $H^1(G,A) \defeq Z^1(G,A) / \sim$. Clearly, $H^1(G,A)$ is (covariantly) functorial in $A$.

\bprop
Let $G \curvearrowright X$ be a topological dynamical system on a compact space $X$. Let $H$ be a group.

$G \curvearrowright X$ is continuous $H$-cocycle rigid if and only if the canonical map $H \to C(X,H)$ (the map dual to $X \to \gekl{\rm pt}$) induces a surjective map $H^1(G,H) \to H^1(G,C(X,H))$.
\eprop
Note that we equip $H$ with the trivial $G$-action, and $G$ acts on $C(X,H)$ via $(s.a)(x) = a(s^{-1}.x)$.
\bproof
Just check that $c: \: C(G,C(X,H)) \to C(G \times X,H)$ defined by $c(a)(g,x) = a_g(g,x)$ is a bijection, with inverse given by $c^{-1}(b)_s(x) = b(s,s^{-1}.x)$. $c$ identifies $Z^1(G,C(X,H))$ with the set of continuous $H$-cocycles for $G \curvearrowright X$ in the sense of Definition~\ref{def-cHc}. In addition, $a \sim a'$ if and only if $c(a) \sim c(a')$ in the sense of Definition~\ref{def-cohom}. It is then easy to see that for $a \in Z^1(G,C(X,H))$, the class $[a] \in H^1(G,C(X,H))$ lies in the image of the canonical map $H^1(G,H) \to H^1(G,C(X,H))$ if and only if $c(a) \sim \rho$ for some group homomorphism $\rho: \: G \to H$.
\eproof

Using the language of non-abelian group cohomology, we now prove a positive result in continuous cocycle rigidity.
\bdefin
A topological dynamical system $G \curvearrowright X$ on a compact space $X$ is called almost $\Zz G$-projective if $C(X,\Zz) \cong \Zz \oplus P$ as $\Zz G$-modules, where the copy of $\Zz$ is given by the constant functions on $X$ and $P$ is a projective $\Zz G$-module.

We call $G \curvearrowright X$ almost $\Zz G$-free if $P$ can be chosen to be $\Zz G$-free.
\edefin
Clearly, if a system is almost $\Zz G$-free, then it is almost $\Zz G$-projective.
\bremark
It is easy to see that $G \curvearrowright X$ is almost $\Zz G$-free if we can find a $\Zz$-basis $\cB$ for $C(X,\Zz)$ with the following properties:
\begin{itemize}
\item $\cB$ is $G$-invariant,
\item $1_X \in \cB$,
\item $G$ acts freely on $\cB \setminus \gekl{1_X}$.
\end{itemize}
\eremark

Topological Bernoulli actions for torsion-free groups turn out to be almost $\Zz G$-free.
\bexample
\label{ex-Bernoulli}
Let $G$ be a torsion-free group and $X_0$ a compact space. Then the Bernoulli action $G \curvearrowright X_0^G$ is almost $\Zz G$-free. Namely, choose a $\Zz$-basis $\cB_0$ for $C(X_0,\Zz)$ with $1_{X_0} \in \cB_0$. This is always possible, see for instance \cite[Proposition~2.12]{CEL}. Then set
$$
  \cB = \menge{\rukl{\bigotimes_{f \in F} b_f} \otimes 1_{X_0^{G \setminus F}}}
  {F \subseteq G \ {\rm finite}, b_f \in \cB_0}.
$$
$\cB$ is a $\Zz$-basis as $C(X_0^G,\Zz) = \bigcup_{F \subseteq G \ {\rm finite}} \rukl{\bigotimes_{f \in F} C(X_0)} \otimes 1_{X_0^{G \setminus F}}$. Obviously, $\gekl{1_{X_0^G}}$ lies in $\cB$. Moreover, $G$ acts freely on $\cB \setminus \gekl{1_{X_0^G}}$ as $G$ is torsion-free.
\eexample

Building on the previous example, we now show that for torsion-free groups, subshifts of full shifts over finite alphabets whose forbidden words avoid a fixed letter are almost $\Zz G$-free.
\bexample
\label{ex-subshift}
Let $G$ be a torsion-free group, $A = \gekl{0, \dotsc, N}$ a finite alphabet and $G \curvearrowright A^G$ the full shift. Elements in $A^G$ are of the form $x = (x_{\gamma})_{\gamma \in G}$, and $g \in G$ acts by $(g.x)_{\gamma} = x_{g^{-1}\gamma}$. For every $G$-invariant closed subset $X$ of $A^G$ we can find a collection $\gekl{F_i}_{i \in I}$ of non-empty finite subsets of $G$ and $x_i \in A^{F_i}$, $i \in I$, such that
$$X = \menge{x = (x_{\gamma})_{\gamma} \in A^G}{{\rm For} \ {\rm every} \ i \in I \ {\rm and} \ g \in G, \pi_{F_i}(g.x) \neq x_i}.$$
Here $\pi_{F_i}$ is the canonical projection $A^G \onto A^{F_i}$. $\gekl{x_i}_{i \in I}$ are called the forbidden words for $X$.
\setlength{\parindent}{0.5cm} \setlength{\parskip}{0cm}

Now assume that $X$ is a $G$-invariant closed subset whose forbidden words $x_i$ satisfy $x_i \in \gekl{1, \dotsc, N}$, i.e., all the forbidden words avoid a fixed letter ($0$ in our case). If that is the case, then we claim that $G \curvearrowright X$ is almost $\Zz G$-free.

Here is the reason: Obviously, $\cB_0 = \gekl{1_{A}, 1_{\gekl{1}}, \dotsc, 1_{\gekl{N}}}$ is a $\Zz$-basis for $C(A,\Zz)$. Given a finite subset $\emptyset \neq F \subseteq G$ and $x = (x_f)_{f \in F} \in \gekl{1, \dotsc, N}^F$, let $b(F,x) = \bigotimes_{f \in F} 1_{\gekl{x_f}} \otimes 1_{A^{G \setminus F}}$. As we have seen in Example~\ref{ex-Bernoulli}, $\cB = \gekl{1_{A^G}} \cup \menge{b(F,x)}{\emptyset \neq F \subseteq G \ {\rm finite}, x \in \gekl{1, \dotsc, N}^F}$ is a $\Zz$-basis for $C(A^G,\Zz)$.

Consider the subspace $C_0(A^G \setminus X,\Zz) \subseteq C(A^G,\Zz)$ of functions vanishing on $X$. An element $\sum z_{F,x} b(F,x) \in C(A^G,\Zz)$ ($z_{F,x} \in \Zz$) lies in $C_0(A^G \setminus X,\Zz)$ if and only if for every $(F,x)$ with $z_{F,x} \neq 0$, there exists $i \in I$ and $g \in G$ with $F_i \subseteq gF$ and $\pi_{F_i}(g.x) = x_i$. Clearly, if $(F,x)$ satisfies this property, then $b(F,x)$ lies in $C_0(A^G \setminus X,\Zz)$. Conversely, suppose that $\sum z_{F,x} b(F,x)$ lies in $C_0(A^G \setminus X,\Zz)$ but there exists $(\ti{F},\ti{x})$ with $z_{\ti{F},\ti{x}} \neq 0$ such that for all $i \in I$ and $g \in G$, $F_i \nsubseteq g\ti{F}$ or $\pi_{F_i}(g.\ti{x}) \neq x_i$. Among all the $(\ti{F},\ti{x})$ with this property, choose a pair such that $\ti{F}$ is minimal. Define $w \in A^G$ by setting $w_{\gamma} = \ti{x}_{\gamma}$ if $\gamma \in \ti{F}$ and $w_{\gamma} = 0$ otherwise. Then $w \in X$, and by our choice of $w$ and $(\ti{F},\ti{x})$, we have $b(\ti{F},\ti{x})(w) = 1$ and $b(F',x')(w) = 0$ for all $(F',x') \neq (\ti{F},\ti{x})$ with $z_{F',x'} \neq 0$. Hence $\rukl{\sum z_{F,x} b(F,x)}(w) = z_{\ti{F},\ti{x}}$, which contradicts that $\sum z_{F,x} b(F,x)$ vanishes on $X$. This shows that
$$\cB_v \defeq \menge{b(F,x)}{{\rm There} \ {\rm is} \ i \in I \ {\rm and} \ g \in G \ {\rm with} \ F_i \subseteq gF \ {\rm and} \ \pi_{F_i}(g.x) = x_i}$$
is a $\Zz$-basis for $C_0(A^G \setminus X,\Zz)$.

The canonical homomorphisms give rise to the exact sequence $0 \to C_0(A^G \setminus X,\Zz) \to C(A^G,\Zz) \to C(X,\Zz) \to 0$. One way to see this would be to apply K-theory to the exact sequence $0 \to C_0(A^G \setminus X) \to C(A^G) \to C(X) \to 0$. Therefore, the image $\cB_X$ of $\cB \setminus \cB_v$ under the canonical projection $C(A^G,\Zz) \onto C(X,\Zz)$ is a $\Zz$-basis for $C(X,\Zz)$. As $\cB_v$ is clearly $G$-invariant, so is $\cB_X$. Moreover, $1_X \in \cB_X$, and $G$ acts freely on $\cB_X \setminus \gekl{1_X} \cong \cB \setminus (\cB_v \cup \gekl{1_{A^G}})$. Therefore, $G \curvearrowright X$ is almost $\Zz G$-free.
\eexample
\setlength{\parindent}{0cm} \setlength{\parskip}{0.5cm}

The systems in Example~\ref{ex-Bernoulli} and Example~\ref{ex-subshift} are not minimal. Here are examples of minimal topological dynamical systems which are almost $\Zz G$-free.
\bexample
A Denjoy homeomorphism is a homeomorphism $\varphi$ of the circle which has no periodic points and is not conjugate to a rigid rotation (see for instance \cite{PSS}). It turns out that there is a unique closed $\varphi$-invariant subspace $\Sigma$ of the circle which is minimal for $\varphi$. $\Sigma$ is a Cantor set. The restriction of $\varphi$ to $\Sigma$ gives rise to a Cantor minimal system $\Zz \curvearrowright \Sigma$. In the proof of \cite[Lemma~6.1]{PSS}, $K_0(C(\Sigma)) \cong C(\Sigma,\Zz)$ is identified as a $\Zz[\Zz]$-module with $\Zz \oplus F$, where $F$ is $\Zz[\Zz]$-free. Hence $\Zz \overset{\varphi}{\curvearrowright} \Sigma$ is almost $\Zz[\Zz]$-free.
\eexample

It would be interesting to find more examples of almost $\Zz G$-projective systems.

Here is why we are interested in the property \an{almost $\Zz G$-projective}:
\bprop
\label{CCR-abelian}
Let $G \curvearrowright X$ be a topological dynamical system on a compact space $X$, and suppose that $G \curvearrowright X$ is almost $\Zz G$-projective. Let $H$ be an abelian $G$-group. Then the canonical map $H \to C(X,H)$ induces injections $H^i(G,H) \into H^i(G,C(X,H))$ for every $i \geq 0$. Moreover, assume that for every $G$-module $M$, $H^1(G,\Zz G \otimes M) \cong \gekl{0}$. Then $G \curvearrowright X$ is continuous $H$-cocycle rigid.
\eprop
\bproof
By assumption, $C(X,\Zz) \cong \Zz \oplus P$ (as $\Zz G$-modules). Then $C(X,H) \cong C(X,\Zz) \otimes H \cong H \oplus (P \otimes H)$. Therefore, for every $i \geq 0$, $C(G,H) \to C(C,C(X,H))$ induces injective maps $H^i(G,H) \into H^i(G,C(X,H))$ as these maps correspond to the canonical inclusions $H^i(G,H) \into H^i(G,H) \oplus H^i(G,P \otimes H)$ under the identification $H^i(G,C(X,H)) \cong H^i(G,H) \oplus H^i(G,P \otimes H)$. Moreover, it is clear that for $i = 1$, $H^1(G,H) \to H^1(G,C(X,H))$ is surjective if and only if $H^1(G,P \otimes H) \cong \gekl{0}$. But since $P$ is a projective $\Zz G$-module, we can find a $\Zz G$-module $Q$ such that $P \oplus Q$ is a free $\Zz G$-module. By assumption, $H^1(G,(P \oplus Q) \otimes H)$ vanishes, since $(P \oplus Q) \otimes H$ is of the form $\Zz G \otimes M$ for some $G$-module $M$. Hence $H^1(G,P \otimes H) \cong \gekl{0}$.
\eproof

\blemma
\label{CCR-extension}
Suppose that $1 \to H' \overset{\iota}{\lori} H \overset{\pi}{\lori} H'' \to 1$ is an exact sequence of groups, and assume that $H'$ is abelian. Let $G \curvearrowright X$ be a topological dynamical system on a compact space $X$, and suppose that $G \curvearrowright X$ is almost $\Zz G$-projective and that $H^1(G,\Zz G \otimes M) \cong \gekl{0}$ for every $G$-module $M$.

If $G \curvearrowright X$ is continuous $H''$-cocycle rigid, then $G \curvearrowright X$ is continuous $H$-cocycle rigid.
\elemma
\bproof
Write $C' = C(X,H')$, $C = C(X,H)$ and $C'' = C(X,H'')$. Let $i: \: C' \to C$ and $p: \: C \to C''$ be the homomorphisms induced by $\iota$ and $\pi$. We get the following commutative diagram with exact rows:
\bgloz
  \xymatrix{
  1 \ar[r] & H' \ar[d]^{\varphi'} \ar[r]^{\iota} & H \ar[d]^{\varphi} \ar[r]^{\pi} & H'' \ar[d]^{\varphi''} \ar[r] & 1   
  \\
  1 \ar[r] & C' \ar[r]^i & C \ar[r]^p & C'' \ar[r] & 1
  }
\egloz
where $\varphi'$, $\varphi$ and $\varphi''$ are the canonical homomorphisms.

Take $x \in H^1(G,C)$. Since $\varphi''_*$ is surjective, we can find $\xi'' \in H^1(G,H'')$ with $\varphi''_*(\xi'') = p_*(x)$. Let us first prove the following

{\bf Claim:} There exists $\zeta \in H^1(G,H)$ with $\pi_*(\zeta) = \xi''$.

{\bf Proof of the claim:} Let $\lambda''$ be a $1$-cocycle of $G$ in $H''$ representing $\xi''$. Lift $\lambda''$ to a map $\lambda: \: G \to H$ such that $\pi \circ \lambda = \lambda''$. Then, as in \cite[Chapter~I, \S~5.6]{Se1997}, define a $2$-cocycle $\lambda'$ of $G$ in $H'$ by setting $\lambda'_{s,t} \defeq \lambda_s s.\lambda_t \lambda_{st}^{-1}$ and let $\Delta(\lambda'') \defeq [\lambda'] \in H^2(G,H')$. Since $p \circ \varphi \circ \lambda = \varphi'' \circ \pi \circ \lambda = \varphi'' \circ \lambda''$, $\varphi \circ \lambda$ is a lift of $\varphi \circ \lambda$. Moreover, $[\varphi'' \circ \lambda''] = \varphi''_* [\lambda''] = \varphi''_*(\xi'') = p_*(x)$ lies in $\img(p_*)$. Therefore, by \cite[Chapter~I, \S~5.6, Proposition~41]{Se1997}, $\Delta(\varphi'' \circ \lambda'') = 0$ in $H^2(G,C')$. Hence $0 = \Delta(\varphi'' \circ \lambda'') = [\varphi' \circ \lambda'] = \varphi'_* [\lambda'] = \varphi'_*(\Delta(\lambda''))$. As $\varphi'_*$ is injective by Proposition~\ref{CCR-abelian}, we obtain $\Delta(\lambda') = 0$. And thus, again by \cite[Chapter~I, \S~5.6, Proposition~41]{Se1997}, $\xi''$ lies in $\img(\pi_*)$. This proves our claim.

So we can find $\zeta \in H^1(G,H)$ with $\pi_*(\zeta) = \xi''$. Then $p_*(\varphi_*(\zeta)) = \varphi''_*(\pi_*(\zeta)) = \varphi''_*(\xi'') = p_*(x)$. Let $\beta$ be a $1$-cocycle of $G$ in $H$ representing $\zeta$. Let $b = \varphi(\beta)$. Then $\varphi_*(\zeta) = [b]$. Twisting by $b$ gives rise to a commutative diagram
\bgl
\label{CD}
  \xymatrix{
  H^1(G,{}_{\beta} H') \ar[d]^{\rukl{\widetilde{\varphi'}}_*} \ar[r]^{({}_{\beta} \iota)_*} 
  & \ker(({}_{\beta} \pi)_*) \ar[d]^{\ti{\varphi}_*} \ar[r]_{\cong}^{\tau_{\beta}} 
  & \pi_*^{-1}(\pi_*(\zeta)) \ar[d]^{\varphi_*}
  \\
  H^1(G,{}_b C') \ar[r]^{({}_b i)_*} 
  & \ker(({}_b p)_*) \ar[r]_{\cong}^{\tau_b} 
  & p_*^{-1}(p_*(x))
  }
\egl
Twisting and the bijections $\tau_b$ are explained in \cite[Chapter~I, \S~5.3]{Se1997} and at the beginning of \cite[Chapter~I, \S~5.4]{Se1997}. By \cite[Chapter~I, \S~5.5, Proposition~38]{Se1997}, $({}_b i)_*: \: H^1(G,{}_b C) \to \ker(({}_b p)_*)$ is surjective. Moreover, $\rukl{\widetilde{\varphi'}}_*: \: H^1(G,{}_{\beta} H') \to H^1(G,{}_b C')$ is surjective by Proposition~\ref{CCR-abelian}. Hence, by commutativity of the left square in \eqref{CD}, $\ti{\varphi}_*: \: \ker(({}_{\beta} \pi)_*) \to \ker(({}_b p)_*)$ is surjective. Commutativity of the right square in \eqref{CD} implies that $\varphi_*: \: \pi_*^{-1}(\pi_*(\zeta)) \to p_*^{-1}(p_*(x))$ is surjective. In particular, there exists $\xi \in \pi_*^{-1}(\pi_*(\zeta)) \subseteq H^1(G,H)$ with $\varphi_*(\xi) = x$.
\eproof

\bcor
\label{CCR-solvable}
Let $G \curvearrowright X$ be a topological dynamical system on a compact space $X$. Suppose that $G \curvearrowright X$ is almost $\Zz G$-projective and that $H^1(G,\Zz G \otimes M) \cong \gekl{0}$ for every $G$-module $M$.

Then $G \curvearrowright X$ is continuous $H$-cocycle rigid for every solvable group $H$.
\ecor
\bproof
We proceed inductively on the length of a series $\gekl{1} = H_0 \subseteq H_1 \subseteq \dotsc \subseteq H_n = H$ with $H_i \triangleleft H$ for all $1 \leq i \leq n$ and $H_i / H_{i-1}$ abelian for all $1 \leq i \leq n$. The case $n=1$ is taken care of by Proposition~\ref{CCR-abelian}. To go from $n-1$ to $n$, consider the series $\gekl{1} = H_1 / H_1 \subseteq H_2 / H_1 \subseteq \dotsc \subseteq H_n / H_1 = H / H_1$. By induction hypothesis, $G \curvearrowright X$ is continuous $H / H_1$-cocycle rigid. Applying Lemma~\ref{CCR-extension} to $1 \to H_1 \to H \to H / H_1 \to 1$, we obtain that $G \curvearrowright X$ is continuous $H$-cocycle rigid.
\eproof

\bremark
\label{H1-duality}
If $G$ is a duality group in the sense of \cite[Chapter~VIII, \S~10]{Bro} with ${\rm cd}(G) \neq 1$, then $H^1(G,\Zz G \otimes M) \cong \gekl{0}$ for every $G$-module $M$.
\eremark

Clearly, Corollary~\ref{CCR-solvable} and Remark~\ref{H1-duality} imply
\btheooz[Theorem~\ref{CCR}]
Let $G$ be a torsion-free group, let $X$ be a compact space, and suppose that $G \curvearrowright X$ is a topological dynamical system which is almost $\Zz G$-projective. Furthermore, assume that $G$ is a duality group in the sense of \cite[Chapter~VIII, \S~10]{Bro} with $cd(G) \neq 1$. Then $G \curvearrowright X$ is continuous $H$-cocycle rigid for every solvable group $H$.
\etheooz

\section{Continuous orbit couples and topological couplings}
\label{sec-coc-tc}

Let us build the bridge between continuous orbit equivalence and topological couplings. Let $G$ and $H$ be groups.
\bdefin
A topological coupling for $G$ and $H$ consists of a locally compact space $\Omega$ with commuting free and proper left $G$- and right $H$-actions which admit compact open fundamental domains $\bar{X}$ (for the $H$-action) and $\bar{Y}$ (for the $G$-action).

A topological coupling is topologically free if the corresponding action $G \times H \curvearrowright \Omega$ is topologically free.
\edefin
We often write $G \curvearrowright \Omega \curvearrowleft H$ for our topological coupling, or $G \underset{\bar{X}}{\curvearrowright} \Omega \underset{\bar{Y}}{\curvearrowleft} H$ if we want to keep track of the fundamental domains. Here, by a fundamental domain $\bar{X}$ for $\Omega \curvearrowleft H$, we mean a subspace $\bar{X} \subseteq \Omega$ such that the inclusion $\bar{X} \into \Omega$ induces a homeomorphism $\bar{X} \cong \Omega / H$. Since we require $\bar{X}$ to be compact and open, this means that $\bar{X} \times H \to \Omega, \, (x,h) \ma xh$ is a homoemorphism (where $H$ carries the discrete topology). For the sake of brevity, we refer to topologically free topological couplings as topologically free couplings.

Moreover, topological couplings $G \underset{\bar{X}_1}{\curvearrowright} \Omega_1 \underset{\bar{Y}_1}{\curvearrowleft} H$ and $G \underset{\bar{X}_2}{\curvearrowright} \Omega_2 \underset{\bar{Y}_2}{\curvearrowleft} H$ are isomorphic if there exists a $G \times H$-equivariant homeomorphism $\Omega_1 \overset{\cong}{\lori} \Omega_2$ sending $\bar{X}_1$ to $\bar{X}_2$ and $\bar{Y}_1$ to $\bar{Y}_2$.

We now introduce a notion which is similar to, but weaker than continuous orbit equivalence.
\bdefin
\label{couples}
Let $G \curvearrowright X$ and $H \curvearrowright Y$ be topological dynamical systems.

A continuous map $p: \: X \to Y$ is called a continuous orbit map if there exists a continuous map $a: \: G \times X \to H$ such that $p(g.x) = a(g,x).p(x)$.

A continuous orbit couple for $G \curvearrowright X$ and $H \curvearrowright Y$ consists of continuous orbit maps $p: \: X \to Y$ and $q: \: Y \to X$ such that there exist continuous maps $g: \: X \to G$ and $h: \: Y \to H$ such that $qp(x) = g(x).x$ and $pq(y) = h(y).y$ for all $x \in X, y \in Y$.

A continuous orbit couple for $G$ and $H$ consists of topological dynamical systems $G \curvearrowright X$ and $H \curvearrowright Y$ on compact spaces $X$ and $Y$ and a continuous orbit couple $(p,q)$ for $G \curvearrowright X$ and $H \curvearrowright Y$.

We call a continuous orbit couple for $G$ and $H$ topologically free if $G \curvearrowright X$ and $H \curvearrowright Y$ are topologically free.
\edefin
Note that if $(p,q)$ is a continuous orbit couple for $G \curvearrowright X$ and $H \curvearrowright Y$ such that $q = p^{-1}$ (i.e., $g \equiv e$ and $h \equiv e$), then $G \curvearrowright X$ and $H \curvearrowright Y$ are continuously orbit equivalent. In that case, we call $(p,q)$ a continuous orbit equivalence.

Continuous orbit couples $(p_i,q_i)$ for $G \curvearrowright X_i$ and $H \curvearrowright Y_i$, $i=1,2$, are isomorphic if there exist a $G$-equivariant homeomorphism $X_1 \overset{\cong}{\lori} X_2$ and an $H$-equivariant homeomorphism $Y_1 \overset{\cong}{\lori} Y_2$ such that the diagrams
$$
  \xymatrix{
  X_1 \ar[d]^{\cong} \ar[r]^{p_1} & Y_1 \ar[d]^{\cong}
  \\
  X_2 \ar[r]^{p_2} & Y_2
  }
\ \ \ \ \ \ 
  \xymatrix{
  Y_1 \ar[d]^{\cong} \ar[r]^{q_1} & X_1 \ar[d]^{\cong}
  \\
  Y_2 \ar[r]^{q_2} & X_2
  }
$$
commute.

The main goal of this section is to prove the following
\btheo
\label{1-1_abc}
Let $G$ and $H$ be groups. There is a one-to-one correspondence between isomorphism classes of topologically free couplings for $G$ and $H$ and isomorphism classes of topologically free continuous orbit couples for $G$ and $H$, with the following additional properties:
\begin{enumerate}
\item[(a)] topological couplings with $\bar{X} = \bar{Y}$ correspond to continuous orbit equivalences; 
\item[(b)] for topologically free couplings $G \underset{\bar{X}_1}{\curvearrowright} \Omega_1 \underset{\bar{X}_1}{\curvearrowleft} H$ and $G \underset{\bar{X}_2}{\curvearrowright} \Omega_2 \underset{\bar{X}_2}{\curvearrowleft} H$, there exists a $G \times H$-equivariant homeomorphism $\Omega_1 \cong \Omega_2$ (which might or might not preserve the fundamental domains) if and only if the cocycles $a_1$, $b_1$, $a_2$ and $b_2$ of the corresponding continuous orbit equivalence satisfy $a_1 \sim a_2$ and $b_1 \sim b_2$;
\item[(c)] for a topologically free coupling $G \underset{\bar{X}}{\curvearrowright} \Omega \underset{\bar{X}}{\curvearrowleft} H$, there exists an isomorphism $\rho: \: G \cong H$ and a $G \times H$-equivariant map $\Omega \to G$ (where $G$ is equipped with the canonical left $G$-action, and the right $H$-action is given by $\rho$) if and only if the cocycle $a$ (see Definition~\ref{couples}) of the corresponding continuous orbit equivalence satisfies $a \sim \rho$ for an isomorphism $\rho: \: G \cong H$ (the same isomorphism as for the coupling).
\end{enumerate}
\etheo
For the proof of this theorem, we will now present explicit constructions of continuous orbit couples out of topological couplings and vice versa. The constructions are really the topological analogues of those in \cite[\S~3]{Fur} (see also \cite{Sha} and \cite{Sau}).

\subsection{From topological couplings to continuous orbit couples}
\label{tc->coc}

Suppose that we are given a topological coupling $G \underset{\bar{X}}{\curvearrowright} \Omega \underset{\bar{Y}}{\curvearrowleft} H$ for groups $G$ and $H$. We write $G \times \Omega \to \Omega, \, (g,x) \ma gx$ and $\Omega \times H \to \Omega, \, (x,h) \ma xh$ for the left $G$-action and right $H$-action.

Set $X \defeq \bar{X}$ and $Y \defeq \bar{Y}$. Define a map $p: \: X \to Y$ by requiring $Gx \cap Y = \gekl{p(x)}$ for all $x \in X$. The intersection on the left hand side is taken in $\Omega$. Since $Y \subseteq \Omega$ is compact and open, $p$ is continuous. Moreover, by construction, there is a continuous map $\gamma: \: X \to G$ with $p(x) = \gamma(x) x$.

We now define a $G$-action, denoted by $G \times X \to X, \, (g,x) \ma g.x$, as follows: For every $g \in G$ and $x \in X$, there exists a unique $\alpha(g,x) \in H$ such that $gx \in X \alpha(g,x)$. Since $X$ is compact and open, $\alpha: \: G \times X \to H$ is continuous. Set $g.x \defeq gx\alpha(g,x)^{-1}$. It is easy to check that this defines a (left) $G$-action on $X$.

Similarly, we define a continuous map $q: \: Y \to X$ by requiring $X \cap yH = \gekl{q(y)}$ for all $y \in Y$, and let $\eta: \: Y \to H$ be the continuous map satisfying $q(y) = y\eta(y)$. To define an $H$-action on $Y$, let $\beta(y,h) \in G$ be such that $yh \in \beta(y,h) Y$. Again, $\beta: \: Y \times H \to G$ is continuous. Set $h.y \defeq \beta(y,h^{-1})^{-1}yh^{-1}$.

Let us check that $(p,q)$ is a continuous orbit couple for $G$ and $H$. To determine $p(g.x) = p(gx\alpha(g,x)^{-1})$, we need to identify $G g x \alpha(g,x)^{-1} \cap Y$. We have
$$G g x \alpha(g,x)^{-1} \ni \beta(\gamma(x)x,\alpha(g,x)^{-1})^{-1} \gamma(x) x \alpha(g,x)^{-1} \in Y,$$
so $p(g.x) = \beta(\gamma(x)x,\alpha(g,x)^{-1})^{-1} \gamma(x) x \alpha(g,x)^{-1} = \alpha(g,x).(\gamma(x)x) = \alpha(g,x).p(x)$. Similarly, in order to identify $q(h.y) = q(\beta(y,h^{-1})^{-1}yh^{-1})$, we need to determine $X \cap \beta(y,h^{-1})^{-1}yh^{-1}H$. As
$$X \ni \beta(y,h^{-1})^{-1} y \eta(y) \alpha(\beta(y,h^{-1})^{-1},y \eta(y))^{-1} \in \beta(y,h^{-1})^{-1}yh^{-1}H,$$
we conclude that $q(y.h) = \beta(y,h^{-1})^{-1} y \eta(y) \alpha(\beta(y,h^{-1})^{-1},y \eta(y))^{-1} = \beta(y,h^{-1})^{-1}.y\eta(y) = \beta(y,h^{-1})^{-1}.q(y)$. Finally, $qp(x) = q(\gamma(x)x) = \gamma(x)x\alpha(\gamma(x),x)^{-1} = \gamma(x).x$ and $pq(y) = p(y\eta(y)) = \beta(y,\eta(y))^{-1} y \eta(y) = \eta(y)^{-1}.y$. All in all, we see that $(p,q)$ is a continuous orbit couple for $G \curvearrowright X$ and $H \curvearrowright Y$ in the sense of Definition~\ref{couples}, with $a(g,x) = \alpha(g,x)$, $b(h,y) = \beta(y,h^{-1})^{-1}$, $g(x) = \gamma(x)$ and $h(y) = \eta(y)^{-1}$.

Note that our coupling does not need to be topologically free for this construction. However, it is clear that $G \curvearrowright \Omega \curvearrowleft H$ is topologically free (i.e., $G \times H \curvearrowright \Omega$ is topologically free) if and only if $G \curvearrowright X$ and $H \curvearrowright Y$ are topologically free.

\subsection{From continuous orbit couples to topological couplings}
\label{coc->tc}

Let $G \curvearrowright X$ and $H \curvearrowright Y$ be topologically free systems on compact spaces $X$ and $Y$. Assume that $(p,q)$ is a continuous orbit couple for $G \curvearrowright X$ and $H \curvearrowright Y$, and let $a$, $b$, $g$ and $h$ be as in Definition~\ref{couples}. Define commuting left $G$- and right $H$-actions on $X \times H$ by $g(x,h) = (g.x,a(g,x)h)$, $(x,h)h' = (x,hh')$. Furthermore, define commuting left $G$- and right $H$-actions on $G \times Y$ by $g'(g,y) = (g'g,y)$ and $(g,y)h = (gb(h^{-1},y)^{-1},h^{-1}.y)$.

A straightforward computation, using the cocycle identities (see Lemma~\ref{cocycle-id}) for $a$ and $b$, shows that $\Theta: \: X \times H \to G \times Y, \, (x,h) \ma (g(x)^{-1}b(h^{-1},p(x))^{-1},h^{-1}.p(x))$ is a $G$- and $H$-equivariant homeomorphism whose inverse is given by $\Theta^{-1}: \: G \times Y \to X \times H, \, (g,y) \ma (g.q(y),a(g,q(y))h(y))$. Therefore, if we set $\Omega = X \times H$ as a $G \times H$-space and set $\bar{X} = X \times \gekl{e}$, $\bar{Y} = \Theta^{-1}(\gekl{e} \times Y)$, then this yields the desired topologically free coupling $G \underset{\bar{X}}{\curvearrowright} \Omega \underset{\bar{Y}}{\curvearrowleft} H$.

Note that topological freeness of $G \curvearrowright X$ and $H \curvearrowright Y$ ensures that $a$ and $b$ satisfy the cocycle identities (as in Lemma~\ref{cocycle-id}), which are needed in the preceding computations.

\subsection{One-to-one correspondence and consequences}

We can now finish the 
\bproof[Proof of Theorem~\ref{1-1_abc}]
It is straightforward to check that the constructions described in \S~\ref{tc->coc} and \S~\ref{coc->tc} are inverse to each other up to isomorphism. For instance, if we start with a topologically free coupling $G \underset{\bar{X}}{\curvearrowright} \Omega \underset{\bar{Y}}{\curvearrowleft} H$, construct a continuous orbit couple and then again a topological coupling, we end up with a coupling of the form $G \underset{\ti{X}}{\curvearrowright} \ti{\Omega} \underset{\ti{Y}}{\curvearrowleft} H$ where $\ti{\Omega} = X \times H \cong G \times Y$, $\ti{X} = X \times \gekl{e}$ and $\ti{Y} \cong \gekl{e} \times Y$. It is then obvious that $\ti{\Omega} = X \times H \to \Omega, \, (x,h) \ma xh$ is an isomorphism of the couplings $G \underset{\ti{X}}{\curvearrowright} \ti{\Omega} \underset{\ti{Y}}{\curvearrowleft} H$ and $G \underset{\bar{X}}{\curvearrowright} \Omega \underset{\bar{Y}}{\curvearrowleft} H$. Conversely, if we start with a continuous orbit couple $(p,q)$ for topologically free systems $G \curvearrowright X$ and $H \curvearrowright Y$, construct a topological coupling and then again a continuous orbit couple, we end up with a continuous orbit couple $(\ti{p},\ti{q})$ for $G \curvearrowright \ti{X}$ and $H \curvearrowright \ti{Y}$ where $\ti{X} = X \times \gekl{e}$ and $\ti{Y} \cong \gekl{e} \times Y$. The canonical isomorphisms $X \cong X \times \gekl{e}$ and $Y \cong \gekl{e} \times Y$ give rise to an isomorphism between $(p,q)$ and $(\ti{p},\ti{q})$.

The additional properties (a), (b) and (c) are also easy to check.
\eproof
In particular, this proves Theorem~\ref{1-1}.

\bcor
\label{cor_coe->qi}
Assume that there exists a topologically free continuous orbit couple for groups $G$ and $H$. If $G$ is finitely generated, then so is $H$, and $G$ and $H$ are quasi-isometric.

In particular, if topologically free systems $G \curvearrowright X$ and $H \curvearrowright Y$ on compact spaces $X$ and $Y$ are continuously orbit equivalent, and if $G$ is finitely generated, then so is $H$, and $G$ and $H$ are quasi-isometric.
\ecor
\bproof
By Theorem~\ref{1-1_abc}, if there exists a topologically free continuous orbit couple for $G$ and $H$, then there exists a (topologically free) topological coupling for $G$ and $H$. Our claim follows from \cite[Chapter~IX, Exercise~34]{dlH}.
\eproof

\section{Conclusions}
\label{sec-concl}

Now we are ready for the proofs of Theorem~\ref{main1} and Theorem~\ref{main2}.
\btheooz[Theorem~\ref{main1}]
Let $G$ be a torsion-free group, and let $H$ be a finitely generated nilpotent group which is not virtually infinite cyclic. Assume that $G \curvearrowright X$ is a topologically free system on a compact space $X$ such that $G \curvearrowright X$ is almost $\Zz G$-projective. Furthermore, let $H \curvearrowright Y$ be a topologically free system. If $G \curvearrowright X$ and $H \curvearrowright Y$ are continuously orbit equivalent, then they must be conjugate. 
\etheooz
\bproof[Proof of Theorem~\ref{main1}]
Corollary~\ref{cor_coe->qi} implies that $G$ is finitely generated, and that $G$ and $H$ are quasi-isometric. So $G$ is virtually nilpotent but not virtually infinite cyclic because $H$ has these properties (see \cite{Gro}). This means that $G$ contains a finitely generated, torsion-free, nilpotent group as a subgroup of finite index. Therefore, by \cite[Chapter~VIII, Proposition~(10.2)]{Bro}, $G$ has to be a duality group. Since $G$ is not virtually infinite cyclic, we know that ${\rm cd}(G) \neq 1$. Thus $G \curvearrowright X$ is continuous $H$-cocycle rigid by Theorem~\ref{CCR}. Now assume that $G \curvearrowright X \ \sim_{\text{\tiny coe}} \ H \curvearrowright Y$. As $G$ is torsion-free and virtually nilpotent, hence amenable, and because $G \curvearrowright X$ is continuous $H$-cocycle rigid, Theorem~\ref{CCR->COER} implies that $G \curvearrowright X \ \sim_{\text{\tiny conj}} \ H \curvearrowright Y$.
\eproof
Here is our second main result, followed by its proof.
\btheooz[Theorem~\ref{main2}]
Let $G$ be a duality group in the sense of \cite[Chapter~VIII, \S~10]{Bro} which is not infinite cyclic, and let $H$ be a finitely generated solvable group. Assume that $G \curvearrowright X$ is a topologically free system on a compact space $X$ such that $G \curvearrowright X$ is almost $\Zz G$-projective. Furthermore, let $H \curvearrowright Y$ be a topologically free system. If $G \curvearrowright X$ and $H \curvearrowright Y$ are continuously orbit equivalent, then they must be conjugate.
\etheooz
\bproof[Proof of Theorem~\ref{main2}]
By Corollary~\ref{cor_coe->qi}, $G$ is finitely generated and quasi-isometric to $H$.  Therefore, $G$ is amenable (see \cite[Chapter~IV, 50.~Geometric properties]{dlH}). Moreover, as $G$ is a duality group, it is torsion-free. As $G$ is not infinite cyclic, we must have $cd(G) \neq 1$. Hence Theorem~\ref{CCR} implies that $G \curvearrowright X$ is continuous $H$-cocycle rigid. Now assume that $G \curvearrowright X \ \sim_{\text{\tiny coe}} \ H \curvearrowright Y$. Since $G$ is torsion-free and amenable, and because $G \curvearrowright X$ is continuous $H$-cocycle rigid, Theorem~\ref{CCR->COER} implies that $G \curvearrowright X \ \sim_{\text{\tiny conj}} \ H \curvearrowright Y$.
\eproof

\end{document}